\pdfoutput=1
\documentclass{article}
\usepackage[margin=1in]{geometry}
\usepackage[utf8]{inputenc}
\usepackage[english]{babel}
\usepackage{amsmath}
\usepackage{amssymb}
\usepackage{amsthm}
\usepackage{graphicx}
\usepackage{url}
\usepackage{placeins}
\usepackage{hyperref}

\theoremstyle{plain}
\newtheorem{theorem}{Theorem}[section]
\newtheorem{proposition}[theorem]{Proposition}
\newtheorem{lemma}[theorem]{Lemma}
\newtheorem{corollary}[theorem]{Corollary}

\newtheorem*{hypothesisU}{Hypothesis $U(k)$}
\newtheorem{inputfact}{Input}

\theoremstyle{definition}
\newtheorem{definition}[theorem]{Definition}

\theoremstyle{remark}
\newtheorem{remark}[theorem]{Remark}

\numberwithin{equation}{section}
\newif\ifnumberalldisplays
\numberalldisplaystrue
\ifnumberalldisplays
  \DeclareRobustCommand{\[}{\begin{equation}}
  \DeclareRobustCommand{\]}{\end{equation}}
  \makeatletter
  \renewenvironment{align*}{\start@align\@ne\st@rredfalse\m@ne}{\endalign}
  \makeatother
\fi

\newcommand{\ZZ}{\mathbb{Z}}

\newcommand{\eps}{\varepsilon}

\newcommand{\padicInts}[1]{\ZZ_{#1}}
\newcommand{\padicUnits}[1]{\ZZ_{#1}^\times}
\newcommand{\haar}[1]{\mu_{#1}}
\newcommand{\val}[1]{v_{#1}}
\newcommand{\kronChar}{\chi_d}
\newcommand{\kron}[1]{\kronChar(#1)}
\newcommand{\repSet}[1]{E_{#1}}
\newcommand{\repCount}[1]{N_{#1}}
\newcommand{\bConst}[1]{b(#1)}
\newcommand{\locSet}[1]{S_{#1}}
\NewDocumentCommand{\diagSet}{m o}{T_{#1}\IfValueT{#2}{^{(#2)}}}
\newcommand{\shifts}{\mathbf h}
\NewDocumentCommand{\locDen}{m O{\{0\}}}{\delta_{#2}(#1)}
\NewDocumentCommand{\truncDen}{m m O{\shifts}}{\delta_{#3}^{(#2)}(#1)}
\NewDocumentCommand{\locFac}{m O{\shifts}}{\lambda_{#2}(#1)}
\NewDocumentCommand{\truncFac}{m m O{\shifts}}{\lambda_{#3}^{(#2)}(#1)}
\NewDocumentCommand{\singSer}{O{\shifts}}{\mathfrak S_{#1}}
\newcommand{\kCorrSym}[1]{\rho_{#1}}
\NewDocumentCommand{\kCorr}{m m O{\shifts}}{\kCorrSym{#1}(#3;#2)}
\newcommand{\repDensity}[1]{\kCorrSym{1}(#1)}

\title{Poisson Statistics for Values of Binary Quadratic Forms in Short Intervals}
\author{Jack Heimrath}
\date{November 24, 2026}

\begin{document}

\maketitle

\begin{abstract}
We study the number of integers represented by a binary quadratic form in short intervals. Extending the work of
Freiberg, Kurlberg, and Rosenzweig on sums of two squares, we formulate an analogue of the Hardy--Littlewood prime
$k$-tuple conjecture for the integers represented by a genus of primitive positive definite binary quadratic forms of a
fixed negative discriminant. We show that this hypothesis implies that the number of integers represented by a genus,
and by each class in the genus, in a random interval $(n,n+y]$ with $n\le x$, where $y$ is $\tau$ times the average
spacing $\sqrt{\log x}/b(d)$, is asymptotically Poisson distributed with parameter $\tau$. The proof uses the method of
moments, with Bernays' asymptotic for the number of represented integers in place of the prime number theorem, and an
adaptation of Ford's simplified argument for the average of the singular series. We also give numerical evidence for the
hypothesis.
\end{abstract}

\tableofcontents

\section{Introduction}

By the prime number theorem, the number of primes up to $x$ is asymptotic to $x/\log x$, so that an interval $(n,n+y]$
with $n\le x$ and $y=\lambda\log x$ contains $\lambda$ primes on average. In 1976 Gallagher~\cite{gallagher} showed
that, assuming a uniform version of the Hardy--Littlewood prime $k$-tuple conjecture, the number of primes in such an
interval is asymptotically Poisson distributed with parameter $\lambda$, as $n$ ranges over the integers up to $x$. In
other words, at the scale of their average spacing the primes behave like independent random points of density
$1/\log x$. Gallagher's proof is by the method of moments: the $k$-th moment of the number of primes in the interval is
a sum over $k$-tuples of shifts, the hypothesis replaces the count of each prime $k$-tuple by its singular series, and
the heart of the argument is an estimate showing that the singular series equals $1$ on average over $k$-tuples of
shifts in a box. Ford~\cite{ford} later found a much simpler proof of this estimate.

For the sums of two squares, the role of the prime number theorem is played by a theorem of Landau~\cite{landau} from
1908: the number of integers up to $x$ that are sums of two squares is asymptotic to $Kx/\sqrt{\log x}$, where
$K=0.7642\ldots$ is the Landau--Ramanujan constant. These integers thus have average spacing about $\sqrt{\log x}/K$
near $x$, and the natural windows are $(n,n+y]$ with $y=\lambda\sqrt{\log x}/K$. Freiberg, Kurlberg, and
Rosenzweig~\cite{freiberg-kurlberg-rosenzweig}, motivated by the level spacings of toral point scatterers, proved the
analogue of Gallagher's theorem in this setting, assuming an analogous Hardy--Littlewood-type hypothesis for sums of two
squares. Their work built on earlier results on the pair correlation of sums of two squares: Connors and
Keating~\cite{connors-keating} conjectured a formula for it and found numerical evidence for Poisson spacings, and
Smilansky~\cite{smilansky} expressed the singular series for pairs as a product of $p$-adic densities and showed that,
assuming (a reformulation of) their conjecture, the second moment of the number of sums of two squares in short
intervals is consistent with Poisson statistics. Smilansky also went beyond sums of two squares: for each of the forms
$x^2+dy^2$ with $d=1,2,3,4,7$, which have class number one, Smilansky formulated an analogous pair correlation
conjecture and showed that, assuming it, the second moment is again consistent with Poisson
statistics~\cite[\S6]{smilansky}. These results involve only pair correlation and the second moment, whereas the theorem
of Freiberg, Kurlberg, and Rosenzweig establishes convergence of all moments. As in Gallagher's work, the analytic heart
of their argument is an estimate for the average of a singular series, and they point out that Ford's simplification
applies there as well.

In this article we extend the theorem of Freiberg, Kurlberg, and Rosenzweig from sums of two squares to the integers
represented by an arbitrary primitive positive definite binary quadratic form, or more precisely by a class or a genus
of such forms of a fixed discriminant $d<0$. Here the role of the prime number theorem is played by the dissertation of
Bernays~\cite{bernays} from 1912, written under Landau, which extends Landau's theorem to all such forms: for each class
$C$ of discriminant $d$, the number $\repCount{C}(x)$ of integers up to $x$ represented by $C$ satisfies:
\[
\repCount{C}(x)\sim \bConst{d}\,\frac x{\sqrt{\log x}}\qquad(x\to\infty),
\]
for a constant $\bConst{d}>0$ that depends only on $d$, and not on $C$. The precise statement, together with a second
result of Bernays comparing a class with its genus, is given in \S\ref{sec:bernays}. The integers represented by $C$
therefore have average spacing about $\sqrt{\log x}/\bConst{d}$ near $x$, and we take windows of exactly this length:
for fixed $\tau>0$ we choose $y$ so that $y\,\repCount{C}(x)/x\to\tau$, that is, $y\sim\tau\sqrt{\log x}/\bConst{d}$, so
that a window $(n,n+y]$ with $n\le x$ contains $\tau$ represented integers on average in the limit. We will show that,
conditionally on a uniform Hardy--Littlewood-type hypothesis, the number of represented integers in such a window is
Poisson distributed with parameter $\tau$. In particular, for the five forms considered by
Smilansky~\cite[\S6]{smilansky}, this passes from the second moment to the full distribution, under a hypothesis on
correlations of every order.

\subsection{Notation and main results}

Throughout, $d<0$ is a discriminant (not assumed fundamental), $C(d)$ denotes the form class group of order $h(d)$, and
there are $g(d)$ genera, each containing $h(d)/g(d)$ classes; see~\cite[Theorem~3.9 and Corollary~3.14]{cox}
and~\cite{cassels}. For a positive definite form $f$ of discriminant $d$ we write $\repSet{f}:=f(\ZZ^2)$ for the set of
integers it represents. This depends only on the class $C$ of $f$, so we write $\repSet{C}:=\repSet{f}$, and for a genus
$R$ we set $\repSet{R}:=\bigcup_{C\subset R}\repSet{C}$. Their counting functions are:
\[
\repCount{C}(x):=\#(\repSet{C}\cap[1,x]),\qquad \repCount{R}(x):=\#(\repSet{R}\cap[1,x]).
\]
Further, $p$ denotes a rational prime, $\val{p}$ the $p$-adic valuation, $\padicInts{p}$ the $p$-adic integers,
$\haar{p}$ normalized Haar measure on $\padicInts{p}$ ($\haar{p}(\padicInts{p})=1$), and $\kron{p}=\left(\frac
dp\right)$ the Kronecker symbol, so that for $p\nmid d$, $p$ is split if $\kron{p}=1$ and inert if $\kron{p}=-1$.

Our main results are the following, in which Hypothesis~$U(k)$ is stated precisely in \S\ref{sec:k-point-correlations}.

\begin{theorem}[Class-level Poisson distribution, conditional on Hypothesis $U(k)$]\label{thm:class}
Fix a discriminant $d<0$, a class $C\in C(d)$, and $\tau>0$, and assume that Hypothesis~$U(k)$ holds for every $k\ge1$
for the genus containing $C$. Let $x$ and $y$ be real parameters tending to infinity in such a way that
$y\,\repCount{C}(x)/x\to\tau$. Then for every fixed integer $m\ge0$:
\[
\frac1x\,\#\big\{1\le n\le x:\ \repCount{C}(n+y)-\repCount{C}(n)=m\big\}\ \longrightarrow\
e^{-\tau}\,\frac{\tau^m}{m!}\qquad(x\to\infty).
\]
\end{theorem}

\begin{theorem}[Genus-level Poisson distribution, conditional on Hypothesis $U(k)$]\label{thm:genus}
Fix a discriminant $d<0$, a genus $R$ of discriminant $d$, and $\tau>0$, and assume that Hypothesis~$U(k)$ holds for
every $k\ge1$ for $R$. Let $x$ and $y$ be real parameters tending to infinity in such a way that
$y\,\repCount{R}(x)/x\to\tau$. Then for every fixed integer $m\ge0$:
\[
\frac1x\,\#\big\{1\le n\le x:\ \repCount{R}(n+y)-\repCount{R}(n)=m\big\}\ \longrightarrow\
e^{-\tau}\,\frac{\tau^m}{m!}\qquad(x\to\infty).
\]
\end{theorem}

Since $\repCount{R}(x)=\repCount{C}(x)(1+o(1))$ for every class $C\subset R$ (\S\ref{sec:bernays},
Inputs~\ref{inp:asymptotic} and~\ref{inp:equidistribution}), the conditions $y\,\repCount{C}(x)/x\to\tau$ and
$y\,\repCount{R}(x)/x\to\tau$ are equivalent, so Theorems~\ref{thm:class} and~\ref{thm:genus} apply to the same windows
and no separate calibration is needed. For $d=-4$ there is a single class, that of $x^2+y^2$, and $\bConst{-4}=K$; in
this case Theorem~\ref{thm:class} recovers part (b) of~\cite[Theorem~1.2]{freiberg-kurlberg-rosenzweig}, on the number
of values in a window, but not part (a), on nearest-neighbor spacings. Our Hypothesis~$U(k)$ is also stronger than the
one assumed there: Freiberg, Kurlberg, and Rosenzweig only ask that the analogous estimate hold on average over dilates
of bounded convex sets, whereas we ask that it hold uniformly. We do this mainly to simplify the argument, and an
averaged hypothesis should also suffice, since the proof of Theorem~\ref{thm:moments} only uses Hypothesis~$U(k)$ after
summing over the shifts.

\subsection{Outline of the proof and organization}

The idea is to split the work into an analytic part and an arithmetic reduction. Essentially all of the analytic work is
concentrated in Theorem~\ref{thm:genus}: local densities, a singular series, and a moment computation.
Theorem~\ref{thm:class} is then deduced from it in \S\ref{sec:class-from-genus}, via a classical equidistribution
theorem of Bernays~\cite{bernays} showing that a class and its genus represent the same integers with an error
negligible at the scale the Poisson statistics live at. We prove this reduction first, so that the rest of the article
can focus exclusively on Theorem~\ref{thm:genus}.

For the average of the singular series in \S\ref{sec:singular-series} we do not use the more involved arguments of
Gallagher~\cite{gallagher} and of Freiberg, Kurlberg, and Rosenzweig~\cite{freiberg-kurlberg-rosenzweig}. Instead we
adapt Ford's simplified approach~\cite{ford} to binary quadratic forms: we truncate the Euler product of the singular
series at a small multiple of $\log H$, replace the remaining finite product by a periodic function, and average it over
complete periods using the Chinese remainder theorem. The price is a worse error term, a relative saving of only
$1/\log\log H$ (Theorem~\ref{thm:ford-average}), but this is more than enough for the method of moments.

The article is organized as follows. In \S\ref{sec:bernays} we state the two inputs from Bernays' dissertation, which we
use without proof, and in \S\ref{sec:class-from-genus} we deduce Theorem~\ref{thm:class} from Theorem~\ref{thm:genus}.
In \S\ref{sec:k-point-correlations} we state Hypothesis~$U(k)$. We compute the local densities in \S\ref{sec:local-odd}
and \S\ref{sec:local-two}, study the singular series and its average in \S\ref{sec:singular-series}, and prove
Theorem~\ref{thm:genus} by the method of moments in \S\ref{sec:proof-genus}. Finally, \S\ref{sec:numerics} contains
numerical evidence for Hypothesis~$U(k)$.

\subsection{Acknowledgments}
I would like to thank my supervisors, Professor Francesco Cellarosi and Professor Brad Rodgers, for their guidance and
support while working on this article. I would also like to thank Professor M. Ram Murty and Dr Nic Fellini for their
many helpful comments during talks at the Queen's University Number Theory Seminar.

\subsection{Use of Large Language Models}
Claude Sonnet 5.5 was used to transcribe and format this article from a set of handwritten notes and was also used to
help with the numerical analysis in \S\ref{sec:numerics}. Claude Opus 5.5 was used to review the article. Claude Opus
4.8 was used to compute the values of $\locDen{2}$ in \S\ref{sec:local-two}. Its derivation relied on results
from class field theory, with which the author is not familiar, so it is not included here. The derivation of
$\locDen{2}$ presented in this article, and all other proofs, were produced without the use of large language models.

\section{Two inputs from Bernays' dissertation}\label{sec:bernays}

The reduction of Theorem~\ref{thm:class} to Theorem~\ref{thm:genus} requires two classical facts due to
P. Bernays~\cite{bernays}. Their proofs are omitted as they are beyond the scope of the present article; a
detailed treatment can be found in~\cite{bernays}.

\begin{inputfact}[Per-class asymptotic]\label{inp:asymptotic}
For any class $C$ of discriminant $d$,
\[
\repCount{C}(x)\ \sim\ \bConst{d}\,\frac{x}{\sqrt{\log x}}\qquad(x\to\infty),
\]
for a constant $\bConst{d}>0$ depending only on $d$ --- the same constant for every class of that discriminant,
regardless of which genus it lies in.
\end{inputfact}

\begin{inputfact}[Class-genus equidistribution]\label{inp:equidistribution}
Let $h=h(d)$ denote the class number. There is an exponent $\theta=\theta(d)>0$ such that for a class $C$ inside a
genus $R$, the difference $\Delta:=\repSet{R}\setminus \repSet{C}$ (which is also the symmetric difference, since
$\repSet{C}\subseteq \repSet{R}$) satisfies
\[
\#(\Delta\cap[1,x])\ =\ O_d\!\left(\frac{x}{(\log x)^{1/2+\theta}}\right).
\]
Any $\theta$ with $0<\theta<1/h$ and $\theta\le1/4$ is admissible. If $h\le2$, then every class is its own genus and
$\Delta=\emptyset$.
\end{inputfact}

Input~\ref{inp:asymptotic} is the asymptotic for a single class together with the equal-density statement
of~\cite[p.~92]{bernays}. Input~\ref{inp:equidistribution} is not stated in exactly this form in~\cite{bernays}; it
follows from~\cite[p.~102]{bernays} (an integer represented by a genus but not by every class in it is an
\emph{exceptional number} for some fundamental class) and~\cite[\S3, p.~117]{bernays} (for each fundamental class, the
exceptional numbers up to $x$ number $O\big(x/(\sqrt{\log x}\,(\log x)^{\vartheta'/h})\big)$ for every
$\vartheta'$ with $0\le\vartheta'<1$ and $\vartheta'\le h/4$), by taking $\theta=\vartheta'/h$.

Combining the two inputs, we obtain essentially for free that for any genus $R$ of discriminant $d$ and any class
$C\subset R$, $\repCount{R}(x)=\repCount{C}(x)(1+o(1))\sim \bConst{d}\,x/\sqrt{\log x}$. Thus the identical constant
$\bConst{d}$ governs every genus and every class of the fixed discriminant $d$ at once, which is exactly what makes a
separate window calibration unnecessary for Theorem~\ref{thm:class} versus Theorem~\ref{thm:genus}
(\S\ref{sec:class-from-genus}). More precisely, Input~\ref{inp:equidistribution} says that a class and its genus differ
by a set whose counting function is smaller than $\repCount{R}(x)$ by a factor $(\log x)^\theta$, so that it is
invisible at the scale of Input~\ref{inp:asymptotic}.

\section{\texorpdfstring{Theorem~\ref{thm:class} from Theorem~\ref{thm:genus}}
{Class-level from genus-level Poisson distribution}}
\label{sec:class-from-genus}

In this section we deduce Theorem~\ref{thm:class} from Theorem~\ref{thm:genus}, using the two inputs of
\S\ref{sec:bernays}. The idea is to couple the window counts of $\repSet{C}$ and $\repSet{R}$: the two sets differ only
by a very sparse set, so their window counts coincide exactly on an event whose probability tends to $1$.

\begin{proof}[Proof (of Theorem~\ref{thm:class} from Theorem~\ref{thm:genus})]
Fix the genus $R$ and a class $C\subset R$, and fix $\tau>0$. By Input~\ref{inp:asymptotic} (combined with
Input~\ref{inp:equidistribution}, as shown in \S\ref{sec:bernays}), $\repDensity{x}:=\repCount{R}(x)/x\to0$ like
$\bConst{d}/\sqrt{\log x}$, and $\repCount{R}(x)=\repCount{C}(x)(1+o(1))$. Hence the hypothesis
$y\,\repCount{C}(x)/x\to\tau$ of Theorem~\ref{thm:class} is equivalent to:
\[
y\cdot \repDensity{x}\to\tau\qquad(\text{so }y\sim\tau\sqrt{\log x}/\bConst{d}).
\]
This matches the calibration used for Theorem~\ref{thm:genus} in \S\ref{sec:proof-genus}.

Let $\Delta:=\repSet{R}\setminus \repSet{C}$. Then $\repSet{C}\subseteq \repSet{R}$, since a class is contained in its
genus, and $\Delta$ is exactly the set on which the two counts could differ. For $n$ drawn uniformly from the integers
$1\le n\le x$, let $Y_x^R$ and $Y_x^C$ denote the counts of $\repSet{R}$, respectively $\repSet{C}$, in the window
$(n,n+y]$. We split the proof into two claims.

Claim 1: We have $Y_x^C=Y_x^R$ \emph{exactly}, on the event that the window $(n,n+y]$ contains no point of $\Delta$ at
all.

Proof of Claim 1: Since $\repSet{C}\subseteq \repSet{R}$ and the two sets differ only by $\Delta$, we get:
\[
Y_x^R-Y_x^C\ =\ \#\big(\Delta\cap(n,n+y]\big)\ \ge0.
\]
Hence the two counts agree precisely when $\Delta\cap(n,n+y]=\emptyset$. This proves Claim 1.

Claim 2: The event in Claim 1 has probability $\to1$ as $x\to\infty$.

Proof of Claim 2: By Input~\ref{inp:equidistribution}, we have $\#(\Delta\cap[1,2x])=O_d\big(x/(\log
x)^{1/2+\theta}\big)$ for some fixed $\theta=\theta(d)>0$. This range suffices, since $y=o(x)$, so that for $x$ large
every window $(n,n+y]$ with $n\le x$ lies in $[1,2x]$. Only the positivity of $\theta$ will be needed.

Now average over the random starting point $n\in(0,x]$. Writing
$\#(\Delta\cap(n,n+y])=\sum_{m\in\Delta}\mathbf{1}[n<m\le n+y]$ and swapping the order of summation, we find:
\[
\sum_{n=1}^x\#(\Delta\cap(n,n+y])\ =\ \sum_{m\in\Delta}\#\{n\in[1,x]:\ m-y\le n<m\}.
\]
For each fixed $m$, the inner count is at most $y$ regardless of $m$: it counts integers in $[m-y,m-1]\cap[1,x]$, and
since $m$ is an integer, the full interval $[m-y,m-1]$ contains exactly $\lfloor y\rfloor\le y$ integers. Moreover, the
inner count is zero once $m>x+y$. So the right-hand side is at most $y\cdot\#(\Delta\cap[1,x+y])\le
y\cdot\#(\Delta\cap[1,2x])$. Since $n$ is uniform on the $\lfloor x\rfloor$ integers $1\le n\le x$, dividing by
$\lfloor x\rfloor\ge x/2$ we obtain:
\begin{align*}
\mathbb E\big[\#(\Delta\cap(n,n+y])\big]\ &\le\ \frac y{\lfloor x\rfloor}\cdot\#(\Delta\cap[1,2x])
\ \le\ \frac{2y}x\cdot\#(\Delta\cap[1,2x]) \\
&=\ O_d\!\left(\frac{y}{(\log x)^{1/2+\theta}}\right)\ \stackrel{(\ast)}{=}\ O_d\!\left(\frac1{(\log x)^{\theta}}\right)
\ \longrightarrow\ 0,
\end{align*}
where $(\ast)$ holds because $y=O(\sqrt{\log x})$. By Markov's inequality it follows that
$\mathbb P\big(\#(\Delta\cap(n,n+y])\ge1\big)\to0$ as $x\to\infty$. This proves Claim 2.

We now put our two claims to work. By Claim 2, with probability tending to $1$ the window contains no point of
$\Delta$, so by Claim 1 we have $Y_x^C=Y_x^R$ exactly on an event of probability $\to1$. Since
$Y_x^R\to\mathrm{Poisson}(\tau)$ in distribution by Theorem~\ref{thm:genus}, the same must be true of $Y_x^C$. Since
$Y_x^C=\repCount{C}(n+y)-\repCount{C}(n)$ takes integer values, this means precisely that
$\mathbb P(Y_x^C=m)=\frac1{\lfloor x\rfloor}\#\{1\le n\le x:\ \repCount{C}(n+y)-\repCount{C}(n)=m\}\to
e^{-\tau}\tau^m/m!$ for every fixed $m\ge0$. Since $\lfloor x\rfloor/x\to1$, the same limit holds with $\lfloor
x\rfloor$ replaced by $x$. This completes the proof of Theorem~\ref{thm:class}, granted Theorem~\ref{thm:genus}.
\end{proof}

\begin{remark}
This coupling argument is why no separate class level hypothesis is necessary. To see why the coupling is needed,
consider the naive alternative of running the method of moments for $Y_x^C$ directly, comparing it with $Y_x^R$ tuple
by tuple. Write $\kCorrSym{k}^C(\shifts;x)$ and $\kCorrSym{k}^R(\shifts;x)$ for the $k$-point correlations
(Definition~\ref{def:k-point-correlation}) of $\repSet{C}$ and $\repSet{R}$. Since $\repSet{C}\subseteq \repSet{R}$, an
integer $n$ counted by $\kCorrSym{k}^R$ but not by $\kCorrSym{k}^C$ must have $n+h_i\in\Delta$ for some $i$, so the
only bound that Input~\ref{inp:equidistribution} provides is:
\[
0\ \le\ \kCorrSym{k}^R(\shifts;x)-\kCorrSym{k}^C(\shifts;x)\ \le\ \frac1x\sum_{i=1}^k\#\{1\le n\le x:\ n+h_i\in\Delta\}
\ \le\ \frac kx\,\#(\Delta\cap[1,2x])\ \ll_{k,d}\ \frac1{(\log x)^{1/2+\theta}}.
\]
The $k$-th factorial moment is a sum of $k$-point correlations over the $\asymp y^k$ tuples of shifts in the window, so
summing this bound over all of them gives a total error of size:
\[
\frac{y^k}{(\log x)^{1/2+\theta}}\ \asymp_{k,d,\tau}\ (\log x)^{(k-1)/2-\theta},
\]
to be compared with a main term of size $\tau^k$. For $k=1$ this error is $(\log x)^{-\theta}\to0$, which is exactly
the estimate used in the proof of Claim~2. For every $k\ge2$, however, the exponent is $(k-1)/2-\theta\ge1/2-\theta>0$,
since $\theta\le1/4$ by Input~\ref{inp:equidistribution}, so the error bound tends to infinity and the naive comparison
fails already for the second moment, whatever the value of $\tau$. To make it work one would need to know how $\Delta$
correlates with $\repSet{R}$, that is, a version of Hypothesis~$U(k)$ for each class separately. The coupling argument
avoids this altogether: it only uses the first moment of $\#(\Delta\cap(n,n+y])$, and since the event $Y^C_x=Y^R_x$ has
probability tending to $1$, the full Poisson distribution transfers for free.
\end{remark}

With this reduction in hand, we turn to Theorem~\ref{thm:genus}. Henceforth we will always drop the subscript $R$ and
write $E:=\repSet{R}$, $N(x):=\repCount{R}(x)$, and so on, for the fixed genus under discussion.

\section{\texorpdfstring{$k$-point correlations and the singular series}{k-point correlations and the singular series}}
\label{sec:k-point-correlations}

In this section we state Hypothesis~$U(k)$, the input on which Theorem~\ref{thm:genus} is conditional. It asserts that
the $k$-point correlations of $E$ are governed by a singular series built from local densities, in the spirit of the
Hardy--Littlewood prime $k$-tuple conjecture. Throughout this section, unless otherwise stated, $R$ is a fixed genus of
discriminant $d<0$, and $E=\repSet{R}$ and $N(x)=\repCount{R}(x)$ as at the end of \S\ref{sec:class-from-genus}.

To measure how likely it is that several integers $n+h_1,\dots,n+h_k$ are simultaneously represented, we count the
corresponding values of $n$. This motivates the following definition.

\begin{definition}[$k$-point correlation]\label{def:k-point-correlation}
For $k\ge1$ and a finite set $\shifts=\{h_1,\dots,h_k\}\subset\ZZ$ of distinct integers, the $k$-point correlation is:
\[
\kCorr{k}{x}:=\frac1x\sum_{1\le n\le x}\mathbf 1_E(n+h_1)\cdots\mathbf 1_E(n+h_k),\qquad
\repDensity{x}:=\kCorr{1}{x}[\{0\}]=\frac{N(x)}x.
\]
\end{definition}

\begin{remark}
Our terminology follows~\cite{freiberg-kurlberg-rosenzweig}, where the analogous quantity $R_k(\shifts;x)$ is called
the \emph{correlation function}. This terminology comes from statistical mechanics and the study of spectral statistics,
which is the setting of~\cite{freiberg-kurlberg-rosenzweig}: there the sums of two squares arise as the energy levels of
the square billiard.
\end{remark}

For a prime $p$, let $\locSet{p}:=\locSet{p}^{(R)}\subset\padicInts{p}$ denote the set of values represented, over
$\padicInts{p}$, by (any form in) $R$; this is well defined since forms in the same genus are $\padicInts{p}$-equivalent
at every $p$. As the continuous image of the compact set $\padicInts{p}^2$ under a polynomial map, $\locSet{p}$ is
compact, hence closed and Haar-measurable. This motivates the following definition.

\begin{definition}[Local density]\label{def:density}
For a prime $p$ and a finite set $\shifts\subset\ZZ$ of distinct integers, the local density is:
\[
\locDen{p}[\shifts]:=\haar{p}\Big(\bigcap_i\left(\locSet{p}-h_i\right)\Big).
\]
\end{definition}

For $\shifts=\{0\}$ this is $\haar{p}(\locSet{p})$. For instance, if $d=-3$ and $R$ is the genus of $x^2+xy+y^2$, then
$S_3$ contains no element congruent to $2$ modulo $3$, so the translates $S_3$, $S_3-1$, $S_3-2$ have empty intersection
and $\locDen{3}[\{0,1,2\}]=0$. We will compute $\locSet{p}$ exactly, at every prime, in
\S\ref{sec:local-odd}--\S\ref{sec:local-two}; knowing $\locSet{p}$ explicitly then immediately shows how a truncated
approximation to $\locDen{p}[\shifts]$ converges to it.

In \S\ref{sec:singular-series} we will average the singular series over $\shifts$ by approximating it with a finite
Euler product, truncated at primes $p\le Y$, whose factors are periodic in $\shifts$. For this we need a finite-level
version of the local densities. This motivates the following definition.

\begin{definition}[Truncated density]\label{def:truncated}
For $\alpha\ge1$ and a tuple $\shifts=(h_1,\dots,h_k)\in\ZZ^k$, whose entries need not be distinct, the truncated
density is:
\[
\truncDen{p}{\alpha}:=p^{-\alpha}\#\big\{0\le a<p^\alpha:\ a+h_i\in \locSet{p}\bmod p^\alpha\ \ \forall i\big\}.
\]
\end{definition}

Thus $\truncDen{p}{\alpha}$ asks that the reduction of $a+h_i$ modulo $p^\alpha$ lie in the reduction of
$\locSet{p}$ modulo $p^\alpha$, for each $i$ separately. Repeated entries are allowed because in
\S\ref{sec:singular-series} we average over all of $(\ZZ/p^\alpha\ZZ)^k$; when the entries are distinct we identify
$\shifts$ with the underlying set, which is harmless since the definition does not depend on their order. We now
assemble the local densities into the singular series.

\begin{definition}[Singular series]\label{def:singular}
For a finite set $\shifts\subset\ZZ$ of $k$ distinct integers, the singular series is:
\[
\singSer:=\prod_p\locFac{p},\qquad\locFac{p}:=
\frac{\locDen{p}[\shifts]}{\locDen{p}^k}.
\]
\end{definition}

By \S\ref{sec:local-odd}, we have $\locDen{p}=1$ for every odd split prime $p\nmid d$, so $\locFac{p}=1$ there for every
$\shifts$ (\S\ref{sec:local-odd} also shows $\locDen{p}[\shifts]=1$ identically at split primes). Thus the only primes
contributing a factor other than $1$ are the odd inert primes ($\kron{p}=-1$, automatically $p\nmid d$ by definition of
$\kronChar$) together with the (finitely many, for fixed $d$) primes dividing $d$ and $p=2$. Convergence of the infinite
product defining $\singSer$ is established in \S\ref{sec:singular-series}. For instance, for the tuple
$\shifts=\{0,1,2\}$ and $d=-3$ considered above we have $\locDen{3}[\shifts]=0$, so the factor at $p=3$ vanishes and
$\singSer=0$. We call $\shifts$ \emph{inadmissible} if $\singSer=0$. For an inadmissible $\shifts$, no $n$ has
$n+h_1,\dots,n+h_k$ all in $E$.

We are now in a position to state the hypothesis.

\begin{hypothesisU}
Fix $k\ge1$. For every fixed $C>0$:
\[
\max_{\shifts}\Big|\frac{\kCorr{k}{x}}{\repDensity{x}^k}-\singSer\Big|
\ \longrightarrow\ 0\qquad(x\to\infty),
\]
the maximum being over all sets $\shifts\subset(0,C\sqrt{\log x}\,]\cap\ZZ$ of $k$ distinct integers.
\end{hypothesisU}

This is the uniform Hardy--Littlewood-type hypothesis on which Theorem~\ref{thm:genus} is conditional. It asserts that
simultaneous representation of $n+h_1,\dots,n+h_k$ behaves, to leading order, like $k$ independent events each of
probability $\repDensity{x}$, corrected by the local singular series, and that this holds not just for a fixed tuple but
uniformly as the tuple ranges over an interval growing with $x$.

\section{Local densities at odd primes}\label{sec:local-odd}

In this section we compute $\locSet{p}$ and $\locDen{p}$ exactly at every odd prime $p$, together with the bounds on
$\locDen{p}[\shifts]$ needed in \S\ref{sec:singular-series}; apart from two standard facts quoted from~\cite{cox}
and~\cite{serre}, the treatment is self-contained. Fix an odd prime $p$. By Lemma~\ref{lem:leading-coefficient} below,
completing the square identifies $\locSet{p}$, up to a measure-preserving unit rescaling, with the image of the single
quadratic map $Q(X,Y)=X^2-dY^2$. The type of $p$ --- split, inert, or dividing $d$ --- corresponds to whether $Q$,
reduced mod $p$, is isotropic, anisotropic, or degenerate, and this trichotomy is all that is needed to compute
$\locSet{p}$ exactly via one recursive mechanism.

\subsection{Setup: completing the square}\label{subsec:completing-square}

\begin{lemma}[{Leading coefficient; see~\cite[Lemmas~2.3 and~2.25]{cox}}]\label{lem:leading-coefficient}
Every primitive form of discriminant $d$ is properly equivalent to one with leading coefficient prime to any given prime
$p$.
\end{lemma}

Fix such a representative $f=ax^2+bxy+cy^2\in R$, $p\nmid a$. Since $p$ is odd, $2a$ is a unit, and completing the
square gives:
\[
4a\cdot f(x,y)=(2ax+by)^2-dy^2.
\]
Notice that the map $(x,y)\mapsto(X,Y)=(2ax+by,y)$ is a bijection $\padicInts{p}^2\to\padicInts{p}^2$ (invert via
$y=Y,\,x=(X-bY)/(2a)$), so that:
\[
4a\cdot \locSet{p}=\diagSet{p}:=\{X^2-dY^2:\ X,Y\in\padicInts{p}\}.
\]
Since $4a$ is a unit, this rescaling preserves Haar measure \emph{exactly at every finite level}: for any $\alpha$, the
map $n\mapsto4an\bmod p^\alpha$ is a bijection of $\ZZ/p^\alpha\ZZ$, so $\#(A\bmod p^\alpha)=\#(4aA\bmod p^\alpha)$ for
any set $A$, and taking $\alpha\to\infty$ we get $\haar{p}(\locSet{p})=\haar{p}(\diagSet{p})$. It therefore suffices to
compute $\diagSet{p}$.

We will repeatedly need the following classical description of the squares of $p$-adic units, which is where the prime
$2$ first differs from the odd primes.

\begin{lemma}[{Squares of units; see~\cite[Ch.~II, \S3.3, Theorems~3 and~4]{serre}}]\label{lem:squares}
Let $u\in\padicUnits{p}$. If $p$ is odd, then $u$ is a square in $\padicInts{p}$ if and only if it is a square modulo
$p$, so that $(\padicUnits{p})^2$ has index $2$ in $\padicUnits{p}$. If $p=2$, then $u$ is a square if and only if
$u\equiv1\pmod8$, so that $(\padicUnits{2})^2=1+8\padicInts{2}$ has index $4$ in $\padicUnits{2}$.
\end{lemma}

We will not prove this here; a detailed treatment can be found in~\cite[Ch.~II, \S3.3]{serre}. For odd $p$ it follows
directly from Hensel's lemma, and the case $p=2$ will only be needed in \S\ref{sec:local-two}.

\subsection{The nonsingular case}\label{subsec:nonsingular}

\begin{lemma}[Nonsingular resolution]\label{lem:nonsingular}
Let $(X_0,Y_0)\in\padicInts{p}^2$ with $\nabla Q(X_0,Y_0)=(2X_0,-2dY_0)\not\equiv(0,0)\pmod p$, and set
$m:=Q(X_0,Y_0)\in\padicInts{p}$. Then for every $m'\in m+p\padicInts{p}$ the following holds. If $2X_0\not\equiv0\pmod
p$, there is a unique $X\equiv X_0\pmod{p}$ with $Q(X,Y_0)=m'$. Similarly, if $2dY_0\not\equiv0\pmod p$, there is a
unique $Y\equiv Y_0\pmod p$ with $Q(X_0,Y)=m'$.
\end{lemma}

In other words, if $m$ is represented by $Q$ at a nonsingular point, then the whole residue class $m+p\padicInts{p}$ is
contained in $\diagSet{p}$.

The experienced reader will recognize Lemma~\ref{lem:nonsingular} as nothing more than Hensel's
lemma~\cite[Ch.~II, \S2.2]{serre}, applied to the one-variable polynomial $Q(X,Y_0)-m'$ obtained by fixing $Y=Y_0$
(respectively to $Q(X_0,Y)-m'$, fixing $X=X_0$). We nevertheless include a proof, since in this case it is particularly
transparent: we run Newton's method for the quadratic polynomial $g(X)=Q(X,Y_0)$, starting from $X_0$, and observe that,
since $g$ is quadratic, its Taylor expansion is exact, so that the precision doubles at every step.

\begin{proof}
Say $2X_0\not\equiv0\pmod p$ (the case $2dY_0\not\equiv0$ is similar, with $-d$ in place of the leading
coefficient~$1$). Fix $Y=Y_0$ and set $g(X):=Q(X,Y_0)$, so that
$g'(X_0)=2X_0\in\padicUnits{p}$. Given $m'\in m+p\padicInts{p}$, define a sequence by:
\[
X_0\ (\text{given}),\qquad X_{n+1}:=X_n-\frac{g(X_n)-m'}{g'(X_n)}.
\]
Since $g$ is quadratic with leading coefficient $1$, its Taylor expansion around any point terminates exactly at the
quadratic term, with no remainder: $g(X_{n+1})=g(X_n)+g'(X_n)(X_{n+1}-X_n)+(X_{n+1}-X_n)^2$. By definition of $X_{n+1}$,
the first two terms on the right cancel $-(g(X_n)-m')$, leaving the exact identity:
\[
g(X_{n+1})-m'=(X_{n+1}-X_n)^2.
\]

Claim 1: We have $X_n\equiv X_0\pmod p$ and $\beta_n:=\val{p}(g(X_n)-m')\ge2^n$, for every $n\ge0$.

Proof of Claim 1: We argue by induction on $n$. For the base case, $X_0\equiv X_0$ trivially, and
$\beta_0=\val{p}(Q(X_0,Y_0)-m')=\val{p}(m-m')\ge1$ since $m'\in m+p\padicInts{p}$. For the inductive step, assume the
claim at $n$. Then $g'(X_n)=2X_n$ is a unit (as $p$ is odd and $X_n\equiv X_0\pmod p$), so
$\val{p}(X_{n+1}-X_n)=\beta_n-\val{p}(g'(X_n))=\beta_n\ge2^n\ge1$. This gives $X_{n+1}\equiv X_n\equiv X_0\pmod p$.
Moreover, the exact identity above gives $\beta_{n+1}=2\val{p}(X_{n+1}-X_n)=2\beta_n\ge2^{n+1}$, which closes the
induction. This proves Claim 1.

We now put Claim 1 to work. It gives $\beta_n\to\infty$, that is, $\val{p}(X_{n+1}-X_n)\to\infty$, so the sequence is
Cauchy and converges to some $X\in\padicInts{p}$, with $X\equiv X_0\pmod p$ since this holds for every $X_n$. Continuity
then gives:
\[
g(X)-m'=\lim_n(g(X_n)-m')=0.
\]
This proves existence.

Claim 2: The solution $X$ is unique among the elements $\equiv X_0\pmod p$.

Proof of Claim 2: Suppose $X'\equiv X_0\pmod p$ also satisfies $g(X')=m'$. Then:
\[
0=g(X)-g(X')=X^2-X'^2=(X-X')(X+X').
\]
Since $X+X'\equiv2X_0\pmod p$ is a unit (the nonsingularity hypothesis), hence nonzero, it must be the case that
$X=X'$. This proves Claim 2 and completes the proof of the lemma.
\end{proof}

\begin{remark}[Higher precision]\label{rem:nonsingular-precision}
The lemma can be sharpened: if $m'\in m+p^\alpha\padicInts{p}$ for some $\alpha\ge1$, then the same argument, with the
bound $\beta_n\ge\alpha\cdot2^n$ in Claim~1 above, shows that the solution satisfies $X\equiv X_0\pmod{p^\alpha}$
(respectively $Y\equiv Y_0\pmod{p^\alpha}$). In other words, values close to $m$ are represented at points close to
$(X_0,Y_0)$. We will not need this improved precision, since only the representation of the class $m+p\padicInts{p}$ is
used below.
\end{remark}

\subsection{\texorpdfstring{The trichotomy mod $p$}{The trichotomy mod p}}\label{subsec:trichotomy}

The form $Q(X,Y)=X^2-dY^2\pmod p$ is a binary quadratic form over $\mathbb F_p$, and its type is read off directly from
$\kron{p}$:

\begin{itemize}
\item $\kron{p}=1$ (split): $d$ is a nonzero square mod $p$, hence a square in $\padicUnits{p}$ by
Lemma~\ref{lem:squares}: $d=s^2$. Then $Q(X,Y)=(X-sY)(X+sY)=UV$ under the bijection $(X,Y)\mapsto(U,V)=(X-sY,X+sY)$
(invertible over $\padicInts{p}$ since $2,s$ are units).
\item $\kron{p}=-1$ (inert): $d$ is a non-square unit, and $Q$ is anisotropic: $Q(X,Y)\equiv0\pmod p$ forces
$X\equiv Y\equiv0\pmod p$ (else $(X/Y)^2\equiv d$, contradicting that $d$ is a non-square).
\item $\kron{p}=0$ ($p\mid d$): writing $d=p^e\,d_0$ ($e=\val{p}(d)\ge1$, $d_0\in\padicUnits{p}$), we have
$Q(X,Y)\equiv X^2\pmod p$ --- degenerate, rank $1$: the gradient's second component $-2dY$ vanishes mod $p$ identically,
for \emph{every} $Y$, since $p\mid d$.
\end{itemize}

\subsection{Split primes}\label{subsec:split}

\begin{proposition}[Split primes]\label{prop:split}
If $\kron{p}=1$, then $\diagSet{p}=\padicInts{p}$, so $\locDen{p}=1$.
\end{proposition}

\begin{proof}
We have $\diagSet{p}=\{UV:\ U,V\in\padicInts{p}\}$. For any target $m$, take $U=1$ and $V=m$.
\end{proof}

Since $\locSet{p}=\padicInts{p}$ outright (Proposition~\ref{prop:split}), translating by any $h_i$ leaves it unchanged:
$\locSet{p}-h_i=\padicInts{p}$. It follows that $\locDen{p}[\shifts]=1$ identically at split primes, for every
$\shifts$, not merely for $\shifts=\{0\}$. There is nothing left to compute, and split primes never enter the singular
series (\S\ref{sec:k-point-correlations}) at all.

\subsection{Inert primes}\label{subsec:inert}

\begin{proposition}[Inert primes]\label{prop:inert}
If $\kron{p}=-1$, then $\diagSet{p}=\{0\}\cup\{m:\ \val{p}(m)\text{ even}\}$, and $\locDen{p}=p/(p+1)$.
\end{proposition}

\begin{proof}
Since $Q$ is anisotropic mod $p$, $Q(X,Y)\equiv0\pmod p$ forces $X\equiv Y\equiv0\pmod p$. Moreover, every
representative of a target $m\not\equiv0\pmod p$ is automatically nonsingular: if $X\equiv Y\equiv0\pmod p$ then
$Q(X,Y)\equiv0\pmod p$, so a representative of a unit target has $(X,Y)\not\equiv(0,0)\pmod p$, and hence
$\nabla Q=(2X,-2dY)\not\equiv0\pmod p$ since $p\nmid2d$. A nondegenerate binary quadratic form over $\mathbb F_p$
represents every element of $\mathbb F_p$ (a classical counting fact), so for every unit $u$ there is a point
$(X_0,Y_0)$ with $Q(X_0,Y_0)\equiv u\pmod p$, which is nonsingular by what we just showed. Lemma~\ref{lem:nonsingular}
then shows that the whole class $Q(X_0,Y_0)+p\padicInts{p}=u+p\padicInts{p}$ is contained in $\diagSet{p}$. Thus
$\diagSet{p}\supseteq\padicUnits{p}$, with $\val{p}(m)=0$ fully represented.

Now let $\val{p}(m)\ge1$. Only singular representatives are possible (shown above), that is, $X\equiv Y\equiv0\pmod p$.
Writing $X=pX_1$ and $Y=pY_1$ we get:
\[
m=p^2\big(X_1^2-dY_1^2\big)=p^2\,Q(X_1,Y_1).
\]
This is the identical problem, with $X_1,Y_1$ ranging freely over $\padicInts{p}$, and it gives the exact recursive
identity:
\[
\diagSet{p}=\padicUnits{p}\ \cup\ p^2\diagSet{p}.
\]
Unrolling, we find $\diagSet{p}=\padicUnits{p}\cup p^2\padicUnits{p}\cup
p^4\padicUnits{p}\cup\cdots\cup\{0\}=\{0\}\cup\{\val{p}\text{ even}\}$, which proves the set identity. Taking measures,
$\mu:=\haar{p}(\diagSet{p})$ satisfies $\mu=(1-1/p)+p^{-2}\mu$ (using
$\haar{p}(p^2\diagSet{p})=p^{-2}\haar{p}(\diagSet{p})$), so that:
\[
\mu=\frac{1-1/p}{1-1/p^2}=\frac1{1+1/p}=\frac p{p+1}.
\]
This completes the proof.
\end{proof}

\subsection{\texorpdfstring{Primes dividing $d$}{Primes dividing d}}\label{subsec:ramified}

Write $d=p^e\,d_0$, $e=\val{p}(d)\ge1$, $d_0\in\padicUnits{p}$. The idea is again to separate the nonsingular
representatives from the singular ones, which give back the same type of problem with the exponent $e$ lowered by $2$.
This produces a recursion in $e$, whose base cases are $e=0$, already treated in \S\ref{subsec:split} and
\S\ref{subsec:inert}, and $e=1$, which we treat directly.

For notational convenience we keep $d_0$ fixed and introduce, for every integer $e'\ge0$, the set:
\[
\diagSet{p}[e']:=\{X^2-p^{e'}d_0Y^2:\ X,Y\in\padicInts{p}\},
\]
so that $\diagSet{p}=\diagSet{p}[e]$. The set $\diagSet{p}[0]$ is the value set of the form $X^2-d_0Y^2$ with $p\nmid
d_0$, which is described by Proposition~\ref{prop:split} if $\left(\tfrac{d_0}{p}\right)=1$ and by
Proposition~\ref{prop:inert} if $\left(\tfrac{d_0}{p}\right)=-1$.

We start with the nonsingular representatives. Since $p\mid p^{e'}d_0$ for $e'\ge1$, the second component of the
gradient vanishes modulo $p$ for \emph{every} $Y$, so nonsingularity reduces to $X\not\equiv0\pmod p$, with $Y$ entirely
unconstrained.

\begin{lemma}[Nonsingular values]\label{lem:nonsingular-values}
Let $e'\ge1$. The values $X^2-p^{e'}d_0Y^2$ with $X\in\padicUnits{p}$ and $Y\in\padicInts{p}$ are exactly the units that
are squares modulo $p$. In particular, they form a set of measure $\tfrac12(1-\tfrac1p)$.
\end{lemma}

\begin{proof}
If $X$ is a unit, then $\val{p}(p^{e'}d_0Y^2)\ge e'\ge1$, so $X^2-p^{e'}d_0Y^2\equiv X^2\pmod p$ is a unit that is a
square modulo $p$. Conversely, if $u$ is a unit that is a square modulo $p$, then it is a square in $\padicInts{p}$ by
Lemma~\ref{lem:squares}, and $u=(\sqrt u)^2-p^{e'}d_0\cdot0^2$ with $\sqrt u$ a unit. Finally, the units that
are squares modulo $p$ make up $(p-1)/2$ of the $p$ residue classes modulo $p$, which gives the measure.
\end{proof}

Note that the set in Lemma~\ref{lem:nonsingular-values} depends neither on $e'$ nor on $d_0$. We now treat the base case
$e=1$, where the singular representatives can be handled directly. The core idea is that the valuations of the two terms
$X^2$ and $pd_0Y^2$ never coincide, so that no cancellation can occur.

\begin{lemma}[The case $e=1$]\label{lem:ramified-one}
Let $m=p^ju$ with $j\ge0$ and $u\in\padicUnits{p}$. Then $m\in \diagSet{p}[1]$ if and only if either $j$ is even and $u$
is a square modulo $p$, or $j$ is odd and $-u/d_0$ is a square modulo $p$. In particular, whether $m$ lies in
$\diagSet{p}[1]$ depends only on $j$ and on $u\bmod p$, and $\haar{p}(\diagSet{p}[1])=\tfrac12$.
\end{lemma}

\begin{proof}
Write $Q(X,Y):=X^2-pd_0Y^2$. Since $\val{p}(X^2)$ is always even and $\val{p}(pd_0Y^2)=1+2\val{p}(Y)$ is always odd,
these valuations can never coincide, and hence:
\[
\val{p}(Q(X,Y))=\min\big(2\val{p}(X),\,1+2\val{p}(Y)\big)\qquad((X,Y)\ne(0,0)).
\]
Assume first that $j=2s$ is even. A representation $Q(X,Y)=m$ forces $\val{p}(X)=s$ and $\val{p}(Y)\ge s$, so that
$u=X_0^2-pd_0Y_1^2\equiv X_0^2\pmod p$ with $X_0:=X/p^s$ and $Y_1:=Y/p^s$. Hence $u$ is a square modulo $p$. Conversely,
if $u$ is a square modulo $p$, then it is a square in $\padicInts{p}$ by Lemma~\ref{lem:squares}, and $(X,Y)=(p^s\sqrt
u,0)$ represents $m$.

Now assume that $j=2s+1$ is odd. A representation forces $\val{p}(Y)=s$ and $\val{p}(X)\ge s+1$, so that
$u=pX_1^2-d_0Y_0^2\equiv-d_0Y_0^2\pmod p$ with $X_1:=X/p^{s+1}$ and $Y_0:=Y/p^s$, and hence $-u/d_0$ is a square modulo
$p$. Conversely, if $-u/d_0$ is a square modulo $p$, then it is a square in $\padicInts{p}$, again by
Lemma~\ref{lem:squares}, and $(X,Y)=(0,p^s\sqrt{-u/d_0})$ represents $m$.

For the measure, note that for each $j$ the admissible units $u$ form half of the units: the squares modulo $p$ if $j$
is even, and their image under multiplication by the unit $-d_0$ if $j$ is odd. Since the elements of valuation exactly
$j$ form a set of measure $(1-\tfrac1p)p^{-j}$, we obtain:
\[
\haar{p}(\diagSet{p}[1])=\sum_{j\ge0}\tfrac12(1-\tfrac1p)p^{-j}=\tfrac12,
\]
which completes the proof.
\end{proof}

For $e\ge2$ the singular representatives lead back to a problem of the same type.

\begin{lemma}[Recursion]\label{lem:ramified-recursion}
For every $e'\ge2$ we have:
\[
\diagSet{p}[e']=\{u\in\padicUnits{p}:\ u\text{ is a square modulo }p\}\ \cup\ p^2\,\diagSet{p}[e'-2],
\]
and the two sets on the right are disjoint.
\end{lemma}

\begin{proof}
By Lemma~\ref{lem:nonsingular-values}, the values at points with $p\nmid X$ form the first set. If $p\mid X$, write
$X=pX_1$. Then:
\[
X^2-p^{e'}d_0Y^2=p^2\big(X_1^2-p^{e'-2}d_0Y^2\big),
\]
and as $X_1$ and $Y$ range freely over $\padicInts{p}$, these values form exactly $p^2\,\diagSet{p}[e'-2]$. The first
set consists of units and the second of elements divisible by $p^2$, so they are disjoint. This completes the proof.
\end{proof}

Applying Lemma~\ref{lem:ramified-recursion} $\lfloor e/2\rfloor$ times reduces $\diagSet{p}[e]$ to $\diagSet{p}[1]$
(Lemma~\ref{lem:ramified-one}) if $e$ is odd, and to $\diagSet{p}[0]$ (Proposition~\ref{prop:split}
or~\ref{prop:inert}) if $e$ is even. Writing $\mu(e'):=\haar{p}(\diagSet{p}[e'])$, Lemmas~\ref{lem:nonsingular-values}
and~\ref{lem:ramified-recursion} give, for $e'\ge2$:
\[
\mu(e')=\tfrac12(1-\tfrac1p)+p^{-2}\mu(e'-2).
\]
This recursion has fixed point $L=\dfrac p{2(p+1)}$, obtained by solving $L=\tfrac12(1-\tfrac1p)+p^{-2}L$. Since
$\mu(e')-L=p^{-2}(\mu(e'-2)-L)$, iterating gives:
\[
\mu(e)-L=p^{-2\lfloor e/2\rfloor}(\mu(e_0)-L),
\]
where $e_0\in\{0,1\}$ is the parity-matched base exponent. Computing $\mu(e_0)-L$ in each of the three cases yields the
following three propositions.

\begin{proposition}[Odd $e$]\label{prop:ramified-odd}
If $e\ge1$ is odd, then:
\[
\locDen{p}=\frac p{2(p+1)}+\frac1{2(p+1)}\,p^{-(e-1)}.
\]
\end{proposition}

\begin{proposition}[Even $e$, split base: $\left(\tfrac{d_0}{p}\right)=1$]\label{prop:ramified-even-split}
If $e\ge2$ is even and $d_0$ is a square mod $p$, then:
\[
\locDen{p}=\frac p{2(p+1)}+\frac{p+2}{2(p+1)}\,p^{-e}.
\]
\end{proposition}

\begin{proposition}[Even $e$, inert base: $\left(\tfrac{d_0}{p}\right)=-1$]\label{prop:ramified-even-inert}
If $e\ge2$ is even and $d_0$ is not a square mod $p$, then:
\[
\locDen{p}=\frac p{2(p+1)}\big(1+p^{-e}\big).
\]
\end{proposition}

\begin{proof}[Proof (of Propositions~\ref{prop:ramified-odd},~\ref{prop:ramified-even-split},
and~\ref{prop:ramified-even-inert})]
For $e$ odd, Lemma~\ref{lem:ramified-one} gives $\mu(1)-L=\tfrac12-\tfrac p{2(p+1)}=\tfrac1{2(p+1)}$, which gives
Proposition~\ref{prop:ramified-odd} directly from $\mu(e)=L+p^{-(e-1)}(\mu(1)-L)$. For $e$ even, we have
$\mu(0)-L=1-\tfrac p{2(p+1)}=\tfrac{p+2}{2(p+1)}$ in the split case (Proposition~\ref{prop:split}), and
$\mu(0)-L=\tfrac p{p+1}-\tfrac p{2(p+1)}=\tfrac p{2(p+1)}$ in the inert case (Proposition~\ref{prop:inert}). Using
$\mu(e)=L+p^{-e}(\mu(0)-L)$, this gives Propositions~\ref{prop:ramified-even-split}
and~\ref{prop:ramified-even-inert}. This completes the proof.
\end{proof}

\begin{remark}
All three formulas converge to the same limit $p/(2(p+1))$ as $e\to\infty$ --- exactly the density of squares in
$\padicInts{p}$. This is not a coincidence: the recursion's fixed point is forced by the nonsingular piece alone (which
never depends on $d_0$ or the split/inert distinction), and as $e\to\infty$ the recursive piece $p^{-2\lfloor
e/2\rfloor}(\mu(e_0)-L)$ --- the only place the base type enters --- decays to $0$, leaving only the type-independent
fixed point.
\end{remark}

\subsection{Regularity at odd primes}\label{subsec:regularity-odd}

The regularity lemmas below control how fast the truncated densities $\truncDen{p}{\alpha}$ approach
$\locDen{p}[\shifts]$. They rest on a simple structural feature shared by the exact descriptions of $\locSet{p}$ above:
every nonzero element $n$ of $\locSet{p}$ comes with a whole residue class around it, whose modulus is only slightly
larger than $p^{\val{p}(n)}$. At a split or inert prime, for instance, $\locSet{p}$ contains every element of the same
valuation as $n$. This motivates the following definition.

\begin{definition}[$c$-stable set]\label{def:stable}
Let $p$ be a prime and $c\ge1$. A set $A\subseteq\padicInts{p}$ is \emph{$c$-stable} if for every $n\in A\setminus\{0\}$
the whole residue class $n+p^{\val{p}(n)+c}\padicInts{p}$ is contained in $A$. In other words, away from $0$, the set
$A$ is a union of residue classes of the form $n+p^{\val{p}(n)+c}\padicInts{p}$.
\end{definition}

The value $0$ plays no role in this definition. Examples of $1$-stable sets are $\padicInts{p}$, $\{0\}$, and any set
defined by a condition on the valuation alone, such as $\{n:\ \val{p}(n)\text{ even}\}$. A simple non-example is given
by the coset $1+p^s\padicInts{p}$ with $s\ge2$, which is not $1$-stable: it contains $1$ but not the class
$1+p\padicInts{p}$. It is, however, $s$-stable.

\begin{lemma}[Closure properties of $c$-stable sets]\label{lem:stable-closure}
Let $c\ge1$.
\begin{enumerate}
\item[(a)] Any union of $c$-stable subsets of $\padicInts{p}$ is $c$-stable, and a $c$-stable set is $c'$-stable for
every $c'\ge c$.
\item[(b)] If $A$ is $c$-stable, then so is $p^2A$.
\item[(c)] If $A$ is $c$-stable and $\lambda\in\padicInts{p}\setminus\{0\}$, then $\{n\in\padicInts{p}:\ \lambda n\in
A\}$ is $c$-stable.
\item[(d)] If $r\in\padicInts{p}$ and $s\ge1$ satisfy $s\le \val{p}(r)+c$ (where $\val{p}(0)=\infty$), then the coset
$r+p^s\padicInts{p}$ is $c$-stable.
\end{enumerate}
\end{lemma}

\begin{proof}
For (a), let $n\in\bigcup_iA_i$ be nonzero. Then $n\in A_i$ for some $i$, and by $c$-stability of $A_i$ we get
$n+p^{\val{p}(n)+c}\padicInts{p}\subseteq A_i\subseteq\bigcup_iA_i$. For the second assertion, observe that
$n+p^{\val{p}(n)+c'}\padicInts{p}\subseteq n+p^{\val{p}(n)+c}\padicInts{p}$ when $c'\ge c$.

For (b), let $n\in p^2A$ be nonzero, and write $n=p^2m$ with $m\in A\setminus\{0\}$. Then $\val{p}(n)=\val{p}(m)+2$, and
so:
\[
n+p^{\val{p}(n)+c}\padicInts{p}=p^2\bigl(m+p^{\val{p}(m)+c}\padicInts{p}\bigr)\subseteq p^2A,
\]
by $c$-stability of $A$.

For (c), put $B:=\{n\in\padicInts{p}:\ \lambda n\in A\}$ and let $n\in B$ be nonzero, so that $\lambda n\in
A\setminus\{0\}$. If $n'\in n+p^{\val{p}(n)+c}\padicInts{p}$, then $\lambda n'\in\lambda n+p^{\val{p}(\lambda
n)+c}\padicInts{p}\subseteq A$ by $c$-stability of $A$, so that $n'\in B$.

For (d), write $C:=r+p^s\padicInts{p}$ and let $n\in C$ be nonzero. If $\val{p}(r)<s$, then $\val{p}(n)=\val{p}(r)$, so
that $s\le \val{p}(n)+c$ and hence $n+p^{\val{p}(n)+c}\padicInts{p}\subseteq n+p^s\padicInts{p}=C$. If $\val{p}(r)\ge
s$, then $C=p^s\padicInts{p}$ and $\val{p}(n)\ge s$, and every element of $n+p^{\val{p}(n)+c}\padicInts{p}$ has
valuation $\val{p}(n)\ge s$, hence lies in $C$. This completes the proof.
\end{proof}

\begin{lemma}[Stability of $\locSet{p}$ at odd primes]\label{lem:stability-odd}
For every odd prime $p$, the set $\locSet{p}$ is $1$-stable.
\end{lemma}

\begin{proof}
By \S\ref{subsec:completing-square} we have $\locSet{p}=\{m:\ 4a\,m\in \diagSet{p}\}$ with $a$ a unit, so by
Lemma~\ref{lem:stable-closure}(c) it is enough to show that $\diagSet{p}$ is $1$-stable. If $\kron{p}=1$ then
$\diagSet{p}=\padicInts{p}$ (Proposition~\ref{prop:split}), and if $\kron{p}=-1$ then
$\diagSet{p}=\{0\}\cup\{\val{p}\text{ even}\}$ (Proposition~\ref{prop:inert}). In both cases, if $m\in \diagSet{p}$ is
nonzero, then $\diagSet{p}$ contains every element of valuation $\val{p}(m)$, and in particular the class
$m+p^{\val{p}(m)+1}\padicInts{p}$, all of whose elements have valuation $\val{p}(m)$. Hence $\diagSet{p}$ is $1$-stable.

Suppose now that $p\mid d$, write $d=p^ed_0$, and let $\diagSet{p}[e']$ be as in \S\ref{subsec:ramified}, so that
$\diagSet{p}=\diagSet{p}[e]$. We show by induction on $e'$ that $\diagSet{p}[e']$ is $1$-stable for every $e'\ge0$. The
base cases are $e'=0$, which is the case just treated (with $d_0$ in place of $d$), and $e'=1$. For $e'\ge2$,
Lemma~\ref{lem:ramified-recursion} writes $\diagSet{p}[e']$ as the union of the units that are squares modulo $p$ and of
$p^2\,\diagSet{p}[e'-2]$. The first set is a union of cosets $r+p\padicInts{p}$ with $r$ a unit, hence $1$-stable by
Lemma~\ref{lem:stable-closure}(d), and the second is $1$-stable by Lemma~\ref{lem:stable-closure}(b) and the induction
hypothesis. By Lemma~\ref{lem:stable-closure}(a) the union is $1$-stable.

It remains to treat $e'=1$. Let $m=p^ju\in \diagSet{p}[1]$ with $u\in\padicUnits{p}$. Every element of
$m+p^{j+1}\padicInts{p}$ has valuation $j$ and unit part congruent to $u$ modulo $p$, so it lies in $\diagSet{p}[1]$ by
Lemma~\ref{lem:ramified-one}. Hence $m+p^{j+1}\padicInts{p}\subseteq \diagSet{p}[1]$, which shows that $\diagSet{p}[1]$
is $1$-stable and completes the proof.
\end{proof}

The next lemma turns stability into a bound on the truncation error for a single index.

\begin{lemma}[Single-index truncation]\label{lem:truncation}
Let $S\subseteq\padicInts{p}$ be a closed $c$-stable set, and put $\delta:=\haar{p}(S)$ and
$\delta^{(\alpha)}:=p^{-\alpha}\,\#(S\bmod p^\alpha)$. Then for every $\alpha\ge1$:
\[
0\ \le\ \delta^{(\alpha)}-\delta\ \le\ \frac{p^{c-1}}{p^\alpha}.
\]
\end{lemma}

\begin{proof}
Fix $\beta\ge1$ and let $a$ be a residue class modulo $p^\beta$ that meets $S$, so that it is counted in
$\delta^{(\beta)}$. Suppose first that $a\not\equiv0\pmod{p^\beta}$ and $\val{p}(a)\le\beta-c$, where $\val{p}(a)<\beta$
denotes the valuation of any representative. Every $m\in a+p^\beta\padicInts{p}$ then satisfies $\val{p}(m)=\val{p}(a)$,
and $a+p^\beta\padicInts{p}\subseteq m+p^{\val{p}(m)+c}\padicInts{p}$. Applying this to an element $m$ of the class that
lies in $S$, we see by $c$-stability that the whole class $a+p^\beta\padicInts{p}$ is contained in $S$. Consequently all
$p$ of its refinements modulo $p^{\beta+1}$ are counted in $\delta^{(\beta+1)}$.

The remaining classes modulo $p^\beta$ are those with $a\equiv0$ or $\val{p}(a)\ge\beta-c+1$, and there are at most
$p^{c-1}$ of them. Each of these classes meets $S$, so at least one of its $p$ refinements modulo $p^{\beta+1}$ is
counted in $\delta^{(\beta+1)}$, and at most $p-1$ refinements are lost, each of measure $p^{-(\beta+1)}$. Therefore:
\[
0\ \le\ \delta^{(\beta)}-\delta^{(\beta+1)}\ \le\ \frac{(p-1)\,p^{c-1}}{p^{\beta+1}}.
\]
Since $S$ is closed, $S=\bigcap_\beta\big(S+p^\beta\padicInts{p}\big)$, and
$\delta^{(\beta)}=\haar{p}(S+p^\beta\padicInts{p})$ decreases to $\delta$. Telescoping over $\beta\ge\alpha$ gives:
\[
\delta^{(\alpha)}-\delta\ \le\ p^{c-1}\sum_{\beta\ge\alpha}\frac{p-1}{p^{\beta+1}}=\frac{p^{c-1}}{p^\alpha},
\]
and the lower bound holds because $\delta^{(\beta)}$ is decreasing in $\beta$. This completes the proof.
\end{proof}

With these two lemmas in hand, the regularity of the multi-index densities follows from a union bound.

\begin{lemma}[Regularity at odd primes]\label{lem:regularity-odd}
For every finite $\shifts\subset\ZZ$ with $\#\shifts=k\ge1$, every odd prime $p$, and every $\alpha\ge1$:
\[
0\ \le\ \truncDen{p}{\alpha}-\locDen{p}[\shifts]\ \le\ \frac{k}{p^\alpha}.
\]
\end{lemma}

The core idea of the proof is to compare the truncated density with the true density through the union bound over the
$k$ shifts, which reduces the problem to the single-index discrepancy controlled by Lemma~\ref{lem:truncation}.

\begin{proof}
Throughout, $\truncDen{p}{\alpha}\ge\locDen{p}[\shifts]$: truncation only ever gains residues, never loses
them, as the construction below makes precise. This gives the lower bound
$0\le\truncDen{p}{\alpha}-\locDen{p}[\shifts]$ directly, so it is enough to prove the upper bound.

For $C\subseteq\padicInts{p}$, write $\widehat C_\alpha$ for its coset-completion at level $\alpha$, the union of
residue classes mod $p^\alpha$ meeting $C$. Write $A:=\bigcap_iB_i$ for $B_i:=\locSet{p}-h_i$, so that
$\locDen{p}[\shifts]=\haar{p}(A)$, and write $\widetilde A_\alpha:=\bigcap_i\widehat{B_i}_\alpha$ for the intersection
of the \emph{individual} completions.

Claim 1: We have $\truncDen{p}{\alpha}=\haar{p}(\widetilde A_\alpha)$.

Proof of Claim 1: By Definition~\ref{def:truncated}, the element $a$ contributes to $\truncDen{p}{\alpha}$
exactly when $a+h_i$ lies in $\locSet{p}\bmod p^\alpha$ for \emph{each} $i$ separately, that is, when
$a\in\widehat{B_i}_\alpha$ for every $i$. This proves Claim 1.

Observe that $\widetilde A_\alpha$ need not equal $\widehat A_\alpha$, the completion of the intersection, which can be
genuinely smaller. For instance, at $\alpha=1$ with several $h_i$ coinciding mod $p$, a residue can pass every
individual check while admitting no single consistent witness, so the distinction is not merely formal.

Claim 2: We have $\widetilde A_\alpha\setminus A\ \subseteq\ \bigcup_{i=1}^k\big(\widehat{B_i}_\alpha\setminus
B_i\big)$.

Proof of Claim 2: Let $a\in\widetilde A_\alpha\setminus A$. Since $a\notin A=\bigcap_iB_i$, it must be the case that
$a\notin B_{i_0}$ for some specific index $i_0$. On the other hand, $a\in\widetilde A_\alpha$ means directly that
$a\in\widehat{B_{i_0}}_\alpha$, by definition of $\widetilde A_\alpha$ as the intersection over every index. Hence
$a\in\widehat{B_{i_0}}_\alpha\setminus B_{i_0}$, which proves Claim 2.

We now put our claims to work. Since $A\subseteq\widetilde A_\alpha$ always (each $B_i\subseteq\widehat{B_i}_\alpha$),
it suffices to bound $\haar{p}(\widetilde A_\alpha\setminus A)$. By Claim 2 and subadditivity of $\haar{p}$ (the union
bound, applied here to a union of complements, not to the $B_i$'s themselves), we get:
\[
\haar{p}(\widetilde A_\alpha\setminus A)\le\sum_{i=1}^k\haar{p}(\widehat{B_i}_\alpha\setminus B_i).
\]
Each $B_i=\locSet{p}-h_i$ is a translate of $\locSet{p}$, so translation invariance of $\haar{p}$ identifies
$\haar{p}(\widehat{B_i}_\alpha\setminus B_i)$ with the single-index discrepancy
$\truncDen{p}{\alpha}[\{0\}]-\locDen{p}$. The set $\locSet{p}$ is closed (it is compact, as noted before
Definition~\ref{def:density}), and it is $1$-stable by Lemma~\ref{lem:stability-odd}.
Lemma~\ref{lem:truncation} with $c=1$ therefore bounds this discrepancy by $1/p^\alpha$, regardless of $h_i$. Summing
the $k$ terms we find:
\[
\truncDen{p}{\alpha}-\locDen{p}[\shifts]=\haar{p}(\widetilde A_\alpha\setminus A)
\ \le\ \sum_{i=1}^k\frac1{p^\alpha}=\frac k{p^\alpha}.
\]
This completes the proof.
\end{proof}

\subsection{\texorpdfstring{Multi-point densities: near $1$ away from $\det\shifts$}{Multi-point densities: near 1
away from det h}}\label{subsec:near-one}

We begin with an elementary bound, which we will use repeatedly, both here and in \S\ref{sec:singular-series}.

\begin{lemma}[Monotonicity]\label{lem:monotonicity}
Let $p$ be a prime, $\alpha\ge1$, and $\shifts$ a tuple of length $k\ge1$. Then
$\truncDen{p}{\alpha}\le\truncDen{p}{\alpha}[\{0\}]$, and therefore
$\truncFac{p}{\alpha}\le\truncDen{p}{\alpha}[\{0\}]^{1-k}$, where:
\[
\truncFac{p}{\alpha}:=\frac{\truncDen{p}{\alpha}}{\truncDen{p}{\alpha}[\{0\}]^k}.
\]
If the entries of $\shifts$ are distinct, then likewise $\locDen{p}[\shifts]\le\locDen{p}$ and
$\locFac{p}\le\locDen{p}^{1-k}$.
\end{lemma}

\begin{proof}
The set $\bigcap_i(\locSet{p}-h_i)$ is contained in $\locSet{p}-h_1$, which has measure $\locDen{p}$ by translation
invariance of $\haar{p}$. Similarly, every $a$ counted in $\truncDen{p}{\alpha}$ satisfies $a+h_1\in \locSet{p}\bmod
p^\alpha$, and $a\mapsto a+h_1$ is a bijection of $\ZZ/p^\alpha\ZZ$. The bounds for $\lambda$ follow on dividing by
$\truncDen{p}{\alpha}[\{0\}]^k$, respectively $\locDen{p}^k$, which completes the proof.
\end{proof}

For the singular series to converge (\S\ref{sec:singular-series}), we need
$\locFac{p}=\locDen{p}[\shifts]/\locDen{p}^k$ to be close to $1$ for all but finitely many primes;
the next two lemmas show this at the odd inert primes.

\begin{lemma}[Generic inert primes]\label{lem:near-one}
Write $\det(\shifts):=\prod_{i<j}(h_j-h_i)$. If $p$ is odd with $\kron{p}=-1$ and $p\nmid\det(\shifts)$ --- that
is, $h_1,\dots,h_k$ are pairwise distinct mod $p$ --- then:
\[
\locDen{p}[\shifts]=1-\frac k{p+1},\qquad\text{so}\qquad\locFac{p}=1+O_k(p^{-2}).
\]
\end{lemma}

\begin{proof}
Since $h_1,\dots,h_k$ are pairwise distinct mod $p$, for any $a$ at most one index $i_0$ can have
$a+h_{i_0}\equiv0\pmod p$. The $k$ events $\{a\equiv-h_i\pmod p\}$ are therefore pairwise disjoint, each of measure
$1/p$. We split into two cases according to whether $a$ lies in one of these events.

\begin{enumerate}
\item[Case 1:] The element $a$ lies off these events, which has measure $1-k/p$. Then every $a+h_i$ is a unit, hence
automatically in $\locSet{p}$ (Proposition~\ref{prop:inert}: $\val{p}=0$ is even). This gives a full contribution of
$(1-k/p)\times1$.

\item[Case 2:] The element $a$ lies on the event $a\equiv-h_{i_0}\pmod p$, which has measure $1/p$, for each of the $k$
choices of $i_0$. Then every $a+h_j$ with $j\ne i_0$ remains a unit (as $h_j\not\equiv h_{i_0}\pmod p$), so those are
automatically fine, and only $a+h_{i_0}$, now known to be $\equiv0\pmod p$, needs checking. Write $a+h_{i_0}=pc$: we
need
$\val{p}(a+h_{i_0})=1+\val{p}(c)$ to be even, that is, $\val{p}(c)$ odd, and by Proposition~\ref{prop:inert} odd
valuations have density exactly $1-p/(p+1)=1/(p+1)$. So this event contributes $\frac1p\cdot\frac1{p+1}$, for each of
the $k$ choices of $i_0$.
\end{enumerate}

Summing the two cases we obtain:
\[
\locDen{p}[\shifts]=\big(1-\tfrac kp\big)+k\cdot\tfrac1p\cdot\tfrac1{p+1}=1-\tfrac kp\Big(1-\tfrac1{p+1}\Big)
=1-\frac kp\cdot\frac p{p+1}=1-\frac k{p+1}.
\]
This completes the proof.
\end{proof}

\begin{lemma}[Upper bound at inert primes]\label{lem:universal-near-one}
If $p$ is odd with $\kron{p}=-1$, then for every $\shifts$ (not just those with $p\nmid\det\shifts$):
\[
\locFac{p}\ \le\ (1+1/p)^{k-1}.
\]
\end{lemma}

\begin{proof}
By Lemma~\ref{lem:monotonicity} and Proposition~\ref{prop:inert},
$\locFac{p}\le\locDen{p}^{1-k}=(1+1/p)^{k-1}$, as required.
\end{proof}

Together, Lemmas~\ref{lem:near-one} and~\ref{lem:universal-near-one} control every odd inert factor of the singular
series: $\locFac{p}=1+O_k(p^{-2})$ at all but the at most $\omega(\det\shifts)$ inert primes dividing
$\det(\shifts)$, where Lemma~\ref{lem:universal-near-one} applies. Since $\locFac{p}=1$ at odd split
primes, only the finitely many primes dividing $2d$ remain. Absolute convergence of the full product
$\singSer=\prod_p\locFac{p}$ also needs $\locFac{p}$ to be bounded above at these primes (it is nonnegative by
definition), which comes from the lower bounds on $\locDen{p}$ in \S\ref{subsec:ramified} and \S\ref{sec:local-two};
these are combined in Corollary~\ref{cor:convergence}.

\section{\texorpdfstring{Local densities at $p=2$}{Local densities at p=2}}\label{sec:local-two}

In this section we compute $\locSet{2}$ and $\locDen{2}$, running parallel to \S\ref{sec:local-odd}; apart from the case
$p=2$ of Lemma~\ref{lem:squares}, quoted from~\cite{serre}, the treatment is self-contained. At $p=2$ the gradient
$(2X,-2dY)$ of $Q(X,Y)=X^2-dY^2$ vanishes mod $2$ \emph{identically}, so there is no analogue of ``nonsingular mod
$p$''. The recursion survives, but is built instead around the parity of $X$ and $Y$, reflecting that squares in
$\padicUnits{2}$ have index $4$, not $2$ (Lemma~\ref{lem:squares}). We will use the description
$(\padicUnits{2})^2=1+8\padicInts{2}$ repeatedly below.

\subsection{Setup: two completing-the-square regimes}\label{subsec:two-setup}

By Lemma~\ref{lem:leading-coefficient}, we may fix $f=ax^2+bxy+cy^2\in R$ with $a$ odd.

Suppose first that $2\nmid d$. Then necessarily $d\equiv1\pmod4$, and $b$ is odd since $d=b^2-4ac$. The substitution
$X=2ax+by,\,Y=y$ satisfies $X\equiv Y\pmod2$ always (as $b$ is odd), and it is only injective on $\padicInts{2}^2$, not
surjective: its image is exactly $\{(X,Y):\ X\equiv Y\pmod2\}$ (given such $(X,Y)$, recover
$y=Y,\,x=(X-bY)/(2a)\in\padicInts{2}$, since the numerator is even and $2a$ has valuation exactly $1$). So we get:
\[
4a\cdot \locSet{2}=\diagSet{2}:=\{X^2-dY^2:\ X\equiv Y\pmod2\}.
\]
Here $4a$ has valuation $2$, so it is not a unit, and this is a genuine $\times4$ measure scaling rather than a
preserving one. Since $a$ itself is already a unit, we find:
\[
\haar{2}(\diagSet{2})=\haar{2}(4a\cdot \locSet{2})=\tfrac14\haar{2}(a\cdot \locSet{2})=\tfrac14\haar{2}(\locSet{2}),
\]
that is, $\haar{2}(\locSet{2})=4\,\haar{2}(\diagSet{2})$.

Suppose now that $4\mid d$ (forced whenever $2\mid d$, since $d\equiv0,1\bmod4$). Write $d=4m$. Then $b=2b'$ is forced,
and $a\cdot f(x,y)=(ax+b'y)^2-my^2$ (expand, using $m=b'^2-ac$) gives a genuine bijection
$\padicInts{2}^2\to\padicInts{2}^2$. There is no parity obstruction this time, since only the unit $a$, not $2a$, is
used, and we obtain:
\[
a\cdot \locSet{2}=\diagSet{2}:=\{X^2-mY^2:\ X,Y\in\padicInts{2}\
\text{unconstrained}\},\qquad\haar{2}(\locSet{2})=\haar{2}(\diagSet{2}).
\]
This is a genuine unit rescaling here, since $a$ itself is odd.

The values of $X^2-m_0Y^2$ with $\val{2}(m_0)\in\{0,1\}$ will be needed repeatedly below, so we record their structure
once. Call $(X,Y)\in\padicInts{2}^2$ primitive if $X,Y$ are not both even.

\begin{lemma}[Primitive values]\label{lem:primitive-values}
Let $m_0\in\padicInts{2}$ with $\val{2}(m_0)\in\{0,1\}$, let $N_0$ be the set of values $X^2-m_0Y^2$ at primitive
$(X,Y)$, and let $N_{oo},N_{oe},N_{eo}$ be the subsets of values at primitive pairs with $(X,Y)$ both odd, $X$ odd and
$Y$ even, and $X$ even and $Y$ odd, respectively. Then:
\begin{enumerate}
\item[(a)] $\{X^2-m_0Y^2:\ X,Y\in\padicInts{2}\}=\{0\}\cup\bigcup_{j\ge0}4^jN_0$, and $N_0=N_{oo}\cup N_{oe}\cup
N_{eo}$.
\item[(b)] $N_{oo}=(1-m_0)+8\padicInts{2}$.
\item[(c)] If $\val{2}(m_0)=0$, then $N_{oe}=1+4\padicInts{2}$ and $N_{eo}=-m_0+4\padicInts{2}$.
\item[(d)] If $\val{2}(m_0)=1$, then $N_{oe}=1+8\padicInts{2}$ and
$N_{eo}=\big(-m_0+16\padicInts{2}\big)\cup\big(4-m_0+16\padicInts{2}\big)$.
\end{enumerate}
\end{lemma}

\begin{proof}
For part (a), a pair $(X,Y)\ne(0,0)$ can be written uniquely as $2^j(X_0,Y_0)$ with $j=\min(\val{2}(X),\val{2}(Y))$ and
$(X_0,Y_0)$ primitive, and then $X^2-m_0Y^2=4^j(X_0^2-m_0Y_0^2)$. Conversely, every $4^jn_0$ with $n_0\in N_0$ arises
this way. The value $0$ arises at $(0,0)$. The second assertion is the definition of the three subsets.

For parts (b)--(d) recall from Lemma~\ref{lem:squares} that the squares of odd elements of $\padicInts{2}$ are exactly
the elements of $1+8\padicInts{2}$.

For part (b), if $X,Y$ are both odd, then $X^2,Y^2$ range independently over $1+8\padicInts{2}$, so $X^2-m_0Y^2$ lies
in:
\[
(1+8\padicInts{2})-m_0(1+8\padicInts{2})=(1-m_0)+8\padicInts{2},
\]
because $8m_0\padicInts{2}\subseteq8\padicInts{2}$. Every element of this coset occurs already with $Y=1$ and $X^2$
varying.

For part (c), let $m_0$ be odd. If $X$ is odd and $Y=2w$ is even, then $X^2-m_0Y^2=X^2-4m_0w^2$, and
$4m_0w^2\equiv0\pmod8$ for $w$ even while $4m_0w^2\equiv4\pmod8$ for $w$ odd. Hence this value lies in
$(1+8\padicInts{2})\cup(5+8\padicInts{2})=1+4\padicInts{2}$, and both cosets occur (take $w=0$, respectively $w=1$ with
$X^2$ varying). If $X=2w$ is even and $Y$ is odd, then $X^2-m_0Y^2=4w^2-m_0Y^2$, where $m_0Y^2$ ranges over
$m_0+8\padicInts{2}$, and $4w^2\equiv0\pmod8$ for $w$ even while $4w^2\equiv4\pmod8$ for $w$ odd. Hence the value lies
in $(-m_0+8\padicInts{2})\cup(4-m_0+8\padicInts{2})=-m_0+4\padicInts{2}$, and both cosets occur (take $w=0$,
respectively $w=1$, with $Y^2$ varying).

For part (d), let $\val{2}(m_0)=1$. If $X$ is odd and $Y=2w$ is even, then $X^2-m_0Y^2=X^2-4m_0w^2$ with
$4m_0w^2\equiv0\pmod8$, so this value lies in $1+8\padicInts{2}$, and all of $1+8\padicInts{2}$ occurs with $w=0$. If
$X=2w$ is even and $Y$ is odd, then $X^2-m_0Y^2=4w^2-m_0Y^2$. Here $Y^2$ ranges over $1+8\padicInts{2}$ by
Lemma~\ref{lem:squares}, so $m_0Y^2$ ranges over $m_0+8m_0\padicInts{2}=m_0+16\padicInts{2}$, and $4w^2\equiv0\pmod{16}$
for $w$ even while $4w^2\equiv4\pmod{16}$ for $w$ odd. Hence the value lies in
$(-m_0+16\padicInts{2})\cup(4-m_0+16\padicInts{2})$, and both cosets occur (take $w=0$, respectively $w=1$, with $Y^2$
varying). This completes the proof.
\end{proof}

\subsection{\texorpdfstring{The case $2\nmid d$}{The case 2 does not divide d}}\label{subsec:two-odd-d}

Suppose first that $d\equiv1\bmod8$ (split). By Lemma~\ref{lem:squares}, $d=s^2$ for a unit $s$, so $Q$ factors as
$X^2-dY^2=(X-sY)(X+sY)$. This substitution has determinant $2s$ (valuation exactly $1$, not a unit), so it is injective
but not surjective on $\padicInts{2}^2$. Tracing parities directly, we see that $X\equiv Y\bmod2$ corresponds exactly to
$U,V$ both even.

\begin{proposition}[Split at $2$]\label{prop:two-split}
$\diagSet{2}=4\padicInts{2}$, so $\haar{2}(\locSet{2})=1$.
\end{proposition}

\begin{proof}
By the parity correspondence above, $\diagSet{2}=\{UV:\ U,V\text{ both even}\}$. Writing $U=2U'$ and $V=2V'$, with
$U',V'$ free since ``both even'' is the only constraint, we get $UV=4U'V'$. As $U',V'$ range unconstrained over
$\padicInts{2}$, the product $U'V'$ hits every target directly (for example, $U'=1,V'=n$ gives $n$), so
$\diagSet{2}=4\padicInts{2}$ exactly. This gives $\haar{2}(\diagSet{2})=1/4$, and by the scaling of
\S\ref{subsec:two-setup} we conclude that $\haar{2}(\locSet{2})=4\cdot\tfrac14=1$.
\end{proof}

Suppose now that $d\equiv5\bmod8$ (inert). Then $d$ is not a square (Lemma~\ref{lem:squares}), so the factorization
trick is unavailable, and we argue by valuation directly.

\begin{proposition}[Inert at $2$]\label{prop:two-inert}
If $d\equiv5\pmod8$, then $\diagSet{2}=\{0\}\cup\{n:\ \val{2}(n)\text{ even},\ \val{2}(n)\ge2\}$, so
$\haar{2}(\locSet{2})=2/3$.
\end{proposition}

\begin{proof}
Let $(X,Y)\ne(0,0)$ with $X\equiv Y\pmod2$, and let $j=\min(\val{2}(X),\val{2}(Y))$. We consider the cases $j=0$ and
$j\ge1$ separately.

If $j=0$, then $X,Y$ are both odd, and by Lemma~\ref{lem:primitive-values}(b) the values $X^2-dY^2$ fill exactly the
coset $(1-d)+8\padicInts{2}=4+8\padicInts{2}$, which is the set of elements of valuation exactly $2$.

If $j\ge1$, then $X,Y$ are both even and the parity condition is automatic: $(X,Y)=2^j(X_0,Y_0)$ with $(X_0,Y_0)$ an
arbitrary primitive pair, and $X^2-dY^2=4^jn_0$ with $n_0\in N_0$, where $N_0$ is the set of
Lemma~\ref{lem:primitive-values} for $m_0=d$. Since $d\equiv5\pmod8$ is odd, parts (b) and (c) of that lemma give:
\[
N_{oo}=4+8\padicInts{2},\qquad N_{oe}=1+4\padicInts{2},\qquad N_{eo}=-d+4\padicInts{2}=3+4\padicInts{2},
\]
so that $N_0=\padicUnits{2}\cup(4+8\padicInts{2})$, the union of the set of elements of valuation $0$ and the set of
elements of valuation $2$. Therefore $4^jN_0$ is the union of the elements of valuation $2j$ and those of valuation
$2j+2$.

Putting the two cases together we find:
\[
\diagSet{2}=\{0\}\ \cup\ \{\val{2}=2\}\ \cup\ \bigcup_{j\ge1}\big(\{\val{2}=2j\}\cup\{\val{2}=2j+2\}\big)
=\{0\}\cup\{n:\ \val{2}(n)\text{ even},\ \val{2}(n)\ge2\}.
\]
The set $\{\val{2}=2j\}$ has measure $\tfrac12\cdot4^{-j}$, so
$\haar{2}(\diagSet{2})=\sum_{j\ge1}\tfrac12\cdot4^{-j}=\tfrac16$ and $\haar{2}(\locSet{2})=4\cdot\tfrac16=\tfrac23$.
This completes the proof.
\end{proof}

Thus, for $2\nmid d$, the prime $2$ behaves like a split or an inert prime in disguise: for $d\equiv1\bmod8$ we have
$\locDen{2}=1$ and $2$ contributes nothing to the singular series, while for $d\equiv5\bmod8$ we have
$\locDen{2}=2/3$, matching the shape $p/(p+1)$ at ``$p=2$'' although it is obtained by mod-$8$ bookkeeping.

\subsection{\texorpdfstring{The case $4\mid d$}{The case 4 divides d}}\label{subsec:two-recursion}

Write $d=4m$. By \S\ref{subsec:two-setup}, $\haar{2}(\locSet{2})=\haar{2}(\diagSet{2})$ for $\diagSet{2}=\{X^2-mY^2:\
X,Y\in\padicInts{2}\}$, unconstrained. Write $m=4^\kappa\,m_0$ with $\kappa\ge0$ and $\val{2}(m_0)\in\{0,1\}$, and set
$\diagSet{2}[\kappa]:=\{X^2-4^\kappa\,m_0 Y^2\}$, so that $\diagSet{2}=\diagSet{2}[\kappa]$.

We first treat the four base types ($\kappa=0$), that is, $\diagSet{2}[0]=\{X^2-m_0Y^2:\ X,Y\in\padicInts{2}\}$ with
$\val{2}(m_0)\in\{0,1\}$. By Lemma~\ref{lem:primitive-values}(a), $\diagSet{2}[0]=\{0\}\cup\bigcup_{j\ge0}4^jN_0$, so it
suffices to determine $N_0$ in each case. We write $\haar{2}(\locSet{2})=\haar{2}(\diagSet{2}[0])$ throughout
(\S\ref{subsec:two-setup}).

\begin{proposition}[Base type $m_0\equiv1\bmod8$]\label{prop:two-base-1}
For $m_0\equiv1\bmod8$: $\diagSet{2}[0]=\padicInts{2}\setminus\{\val{2}=1\}$, $\haar{2}(\locSet{2})=3/4$.
\end{proposition}

\begin{proof}
Since $1-m_0\equiv0\pmod8$, Lemma~\ref{lem:primitive-values}(b) gives $N_{oo}=8\padicInts{2}$, and part (c) gives
$N_{oe}=1+4\padicInts{2}$ and $N_{eo}=-m_0+4\padicInts{2}=3+4\padicInts{2}$. Thus
$N_0=\padicUnits{2}\cup8\padicInts{2}$: all units, together with all elements of valuation at least $3$, and $0$. Then
$j=0$ contributes every unit and every element of valuation $\ge3$, and $j=1$ contributes
$4\padicUnits{2}=\{\val{2}=2\}$. Together these (and $0$) exhaust $\padicInts{2}\setminus\{\val{2}=1\}$. Conversely, an
element of $4^jN_0$ has valuation $2j$ or at least $2j+3$, which is never $1$. Hence
$\diagSet{2}[0]=\padicInts{2}\setminus\{\val{2}=1\}$, and $\haar{2}(\{\val{2}=1\})=\tfrac12\cdot\tfrac12=\tfrac14$ gives
$\haar{2}(\locSet{2})=\tfrac34$.
\end{proof}

\begin{proposition}[Base type $m_0\equiv5\bmod8$]\label{prop:two-base-5}
For $m_0\equiv5\bmod8$: $\diagSet{2}[0]=\{0\}\cup\{n:\ \val{2}(n)\text{ even}\}$, $\haar{2}(\locSet{2})=2/3$.
\end{proposition}

\begin{proof}
Since $1-m_0\equiv4\pmod8$, Lemma~\ref{lem:primitive-values} gives $N_{oo}=4+8\padicInts{2}$, $N_{oe}=1+4\padicInts{2}$,
and $N_{eo}=-m_0+4\padicInts{2}=3+4\padicInts{2}$, so $N_0=\padicUnits{2}\cup(4+8\padicInts{2})$, the elements of
valuation $0$ and those of valuation $2$. As in the proof of Proposition~\ref{prop:two-inert}, $4^jN_0$ consists of the
elements of valuation $2j$ and $2j+2$, so $\diagSet{2}[0]=\{0\}\cup\{\val{2}\text{ even}\}$ and
$\haar{2}(\diagSet{2}[0])=\sum_{j\ge0}\tfrac12\cdot4^{-j}=\tfrac23$.
\end{proof}

\begin{proposition}[Base type $m_0\equiv3\bmod4$]\label{prop:two-base-37}
For $m_0\equiv3,7\bmod8$ (i.e.\ $m_0$ odd, $m_0\equiv3\bmod4$):
\[
\diagSet{2}[0]\setminus\{0\}=\bigcup_{j\ge0}4^j\Big((1+4\padicInts{2})\cup\big((1-m_0)+8\padicInts{2}\big)\Big),
\]
where the second coset consists of elements of valuation exactly $1$, and $\haar{2}(\locSet{2})=1/2$.
\end{proposition}

\begin{proof}
Here $-m_0\equiv1\pmod4$, so Lemma~\ref{lem:primitive-values}(c) gives $N_{oe}=N_{eo}=1+4\padicInts{2}$, and part (b)
gives $N_{oo}=(1-m_0)+8\padicInts{2}$, a single coset mod $8$ of elements of valuation exactly $1$ (as
$1-m_0\equiv2\pmod4$). So $N_0$ is the disjoint union of $1+4\padicInts{2}$, of measure $\tfrac14$ and consisting of
units, and $(1-m_0)+8\padicInts{2}$, of measure $\tfrac18$ and consisting of elements of valuation $1$. By
Lemma~\ref{lem:primitive-values}(a), this gives the stated description of $\diagSet{2}[0]$. Moreover $4^jN_0$ consists
of elements of valuations $2j$ and $2j+1$ only, the sets $4^jN_0$ for different $j$ are pairwise disjoint, and
$\haar{2}(4^jN_0)=4^{-j}\cdot\tfrac38$. Therefore
$\haar{2}(\diagSet{2}[0])=\sum_{j\ge0}4^{-j}\cdot\tfrac38=\tfrac38\cdot\tfrac43=\tfrac12$.
\end{proof}

\begin{proposition}[Base type $\val{2}(m_0)=1$]\label{prop:two-base-even}
For $\val{2}(m_0)=1$:
\[
\diagSet{2}[0]\setminus\{0\}=\bigcup_{j\ge0}4^j\Big(\big(\{1,1-m_0\}+8\padicInts{2}\big)
\cup\big(\{-m_0,4-m_0\}+16\padicInts{2}\big)\Big),
\]
where the first two cosets consist of units and the last two of elements of valuation exactly $1$, and
$\haar{2}(\locSet{2})=1/2$.
\end{proposition}

\begin{proof}
By Lemma~\ref{lem:primitive-values}(b) and (d), $N_0$ is the union of the two sets
$(1+8\padicInts{2})\cup\big((1-m_0)+8\padicInts{2}\big)$ and $(-m_0+16\padicInts{2})\cup(4-m_0+16\padicInts{2})$, which
gives the stated description of $\diagSet{2}[0]$ by part (a) of that lemma. The first set consists of units, and its two
cosets are distinct because $1-m_0\equiv3\pmod4$; it has measure $\tfrac18+\tfrac18=\tfrac14$. The second consists of
elements of valuation exactly $1$, and its two cosets mod $16$ are distinct (they differ by $4$); it has measure
$\tfrac1{16}+\tfrac1{16}=\tfrac18$. As in the previous proof, the sets $4^jN_0$ are pairwise disjoint, with
$\haar{2}(4^jN_0)=4^{-j}\cdot\tfrac38$, so $\haar{2}(\diagSet{2}[0])=\tfrac38\cdot\tfrac43=\tfrac12$.
\end{proof}

With the four base types in hand, we now turn to $\kappa\ge1$, where we split on the parity of $X$ alone. If $X$ is
even, then $X=2X_1$, and we find:
\[
X^2-4^\kappa m_0Y^2=4\big(X_1^2-4^{\kappa-1}m_0Y^2\big),
\]
which is the same problem one level down in $\kappa$, with $X_1,Y$ free. This motivates the following description.

\begin{proposition}[The recursion at $2$]\label{prop:two-recursion-sets}
Let $\kappa\ge1$, and let $B_\kappa$ be the set of values of $X^2-4^\kappa m_0Y^2$ at pairs $(X,Y)$ with $X$ odd. Then:
\[
\diagSet{2}[\kappa]=B_\kappa\cup4\,\diagSet{2}[\kappa-1],\qquad
B_\kappa=\begin{cases}
1+4\padicInts{2}&\text{if }\kappa=1\text{ and }m_0\text{ is odd},\\
1+8\padicInts{2}&\text{otherwise},
\end{cases}
\]
and the two sets in the union are disjoint.
\end{proposition}

\begin{proof}
The values at pairs with $X$ even form $4\,\diagSet{2}[\kappa-1]$, by the computation above. If $X$ is odd (any $Y$),
then $\val{2}(4^\kappa\,m_0 Y^2)=2\kappa+\val{2}(m_0)+2\val{2}(Y)\ge2\kappa\ge2$, while $\val{2}(X^2)=0$. These can
never coincide, so every such value is a unit. For $\kappa\ge2$, we have $4^\kappa\,m_0 Y^2\equiv0\pmod8$ regardless of
$Y$ (as $2\kappa\ge4$), so $X^2-4^\kappa m_0Y^2\equiv X^2\equiv1\pmod8$, and every element of $1+8\padicInts{2}$ is
attained already with $Y=0$, since it is a square by Lemma~\ref{lem:squares}. For $\kappa=1$, an identical computation
splits on the parity of $Y$, giving $1+8\padicInts{2}$ if $Y$ is even and $(1-4m_0)+8\padicInts{2}$ if $Y$ is odd. For
$m_0$ odd we have $4m_0\equiv4\pmod8$, so the union is $(1+8\padicInts{2})\cup(5+8\padicInts{2})=1+4\padicInts{2}$,
while for $\val{2}(m_0)=1$ we have $4m_0\equiv0\pmod8$, so both parities give $1+8\padicInts{2}$. Finally, $B_\kappa$
consists of units and $4\,\diagSet{2}[\kappa-1]$ of elements divisible by $4$, so the two sets are disjoint. This
completes the proof.
\end{proof}

\begin{proposition}[The local density at $2$ when $4\mid d$]\label{prop:two-recursion}
Write $\nu(\kappa):=\haar{2}(\diagSet{2}[\kappa])$. By Proposition~\ref{prop:two-recursion-sets}, the recursion
$\nu(\kappa)=\haar{2}(B_\kappa)+\tfrac14\nu(\kappa-1)$ holds for $\kappa\ge1$, with $\haar{2}(B_\kappa)=1/8$ for
$\kappa\ge2$ and $\haar{2}(B_1)=1/4$ (odd $m_0$) or $1/8$ (even $m_0$). Solving: the fixed point of the $\kappa\ge2$
recursion is $L=1/6$, and:
\[
\nu(\kappa)=\frac16+\beta\cdot4^{-\kappa},\qquad\beta=4\Big(\nu(1)-\frac16\Big),\qquad \kappa\ge1,
\]
with $\nu(1)=\haar{2}(B_1)+\tfrac14\nu(0)$ computed from Propositions~\ref{prop:two-base-1}--\ref{prop:two-base-even},
giving:
\[
\beta=\begin{cases}
13/12&m_0\equiv1\bmod8\ (\text{split base})\\
1&m_0\equiv5\bmod8\ (\text{inert base})\\
5/6&m_0\equiv3,7\bmod8\ (\text{ramified-odd base})\\
1/3&\val{2}(m_0)=1\ (\text{even-}m_0\text{ base})
\end{cases}.
\]
\end{proposition}

\begin{proof}
We first treat $\kappa\ge2$. Substituting $\haar{2}(B_\kappa)=1/8=\tfrac16\cdot\tfrac34$ and solving, we get
$\nu(\kappa)-\tfrac16=\tfrac14\big(\nu(\kappa-1)-\tfrac16\big)$, so that
$\nu(\kappa)-\tfrac16=4^{-(\kappa-1)}(\nu(1)-\tfrac16)$, that is, $\nu(\kappa)=\tfrac16+\beta\cdot4^{-\kappa}$ with
$\beta=4(\nu(1)-\tfrac16)$, for every $\kappa\ge1$ (the formula for $\kappa=1$ recovers itself trivially). It remains to
compute $\nu(1)$ in each of the four cases.

\begin{enumerate}
\item[Case 1:] Split base. We have $\nu(1)=\tfrac14+\tfrac14\cdot\tfrac34=\tfrac7{16}$, hence
$\beta=4(\tfrac7{16}-\tfrac16)=4\cdot\tfrac{13}{48}=\tfrac{13}{12}$.
\item[Case 2:] Inert base. We have $\nu(1)=\tfrac14+\tfrac14\cdot\tfrac23=\tfrac5{12}$, hence
$\beta=4(\tfrac5{12}-\tfrac16)=4\cdot\tfrac14=1$.
\item[Case 3:] Ramified-odd base. We have $\nu(1)=\tfrac14+\tfrac14\cdot\tfrac12=\tfrac38$, hence
$\beta=4(\tfrac38-\tfrac16)=4\cdot\tfrac5{24}=\tfrac56$.
\item[Case 4:] Even-$m_0$ base. We have $\nu(1)=\tfrac18+\tfrac14\cdot\tfrac12=\tfrac14$, hence
$\beta=4(\tfrac14-\tfrac16)=4\cdot\tfrac1{12}=\tfrac13$.
\end{enumerate}

This completes the proof.
\end{proof}

Since $\haar{2}(\locSet{2})=\haar{2}(\diagSet{2})$ in this regime (\S\ref{subsec:two-setup}, genuine unit rescaling),
Proposition~\ref{prop:two-recursion} gives $\locDen{2}=\nu(\kappa)=\tfrac16+\beta\cdot4^{-\kappa}$ directly, for
every conductor depth $\kappa\ge1$. This completes the exact computation of $\locDen{2}$ in every case.

\subsection{\texorpdfstring{Regularity at $p=2$}{Regularity at p=2}}\label{subsec:regularity-two}

At $p=2$ the set $\locSet{2}$ need no longer contain, together with each nonzero element $m$, the whole class
$m+2^{\val{2}(m)+1}\padicInts{2}$, since squares of units are detected modulo $8$ and not modulo $2$. It does, however,
contain the class $m+2^{\val{2}(m)+3}\padicInts{2}$; that is, $\locSet{2}$ is $3$-stable in the sense of
Definition~\ref{def:stable}.

\begin{lemma}[Stability of $\locSet{2}$]\label{lem:stability-two}
The set $\locSet{2}$ is $3$-stable, and it is $1$-stable if $d$ is odd.
\end{lemma}

\begin{proof}
By \S\ref{subsec:two-setup} we have $\locSet{2}=\{m:\ \lambda m\in \diagSet{2}\}$ with $\lambda=4a$ if $2\nmid d$ and
$\lambda=a$ if $4\mid d$, so by Lemma~\ref{lem:stable-closure}(c) it is enough to show that $\diagSet{2}$ is $3$-stable,
and $1$-stable if $d$ is odd.

Suppose first that $2\nmid d$. By Propositions~\ref{prop:two-split} and~\ref{prop:two-inert}, $\diagSet{2}$ is either
$4\padicInts{2}$ or $\{0\}\cup\{n:\ \val{2}(n)\text{ even},\ \val{2}(n)\ge2\}$. In both cases, if $n\in \diagSet{2}$ is
nonzero, then $\diagSet{2}$ contains every element of valuation $\val{2}(n)$, and in particular the class
$n+2^{\val{2}(n)+1}\padicInts{2}$. Hence $\diagSet{2}$ is $1$-stable, and so also $3$-stable by
Lemma~\ref{lem:stable-closure}(a).

Suppose now that $4\mid d$, and write $\diagSet{2}=\diagSet{2}[\kappa]$ as in \S\ref{subsec:two-recursion}. We argue by
induction on $\kappa$. Let $\kappa=0$. By Lemma~\ref{lem:primitive-values}(a),
$\diagSet{2}[0]=\{0\}\cup\bigcup_{j\ge0}4^jN_0$ and $N_0=N_{oo}\cup N_{oe}\cup N_{eo}$. By parts (b)--(d) of that lemma,
each of $N_{oo}$, $N_{oe}$, and $N_{eo}$ is a union of cosets $r+2^s\padicInts{2}$ of one of the following kinds:
\begin{itemize}
\item cosets of $8\padicInts{2}$ or $4\padicInts{2}$ with $r$ a unit, so that $s\le3=\val{2}(r)+3$;
\item the coset $(1-m_0)+8\padicInts{2}$ when $\val{2}(m_0)=0$, where $\val{2}(1-m_0)\ge1$, so that $s=3\le
\val{2}(1-m_0)+3$;
\item the cosets $-m_0+16\padicInts{2}$ and $4-m_0+16\padicInts{2}$ when $\val{2}(m_0)=1$, both of whose representatives
have valuation exactly $1$, so that $s=4=1+3$.
\end{itemize}
Each of these is $3$-stable by Lemma~\ref{lem:stable-closure}(d), so $N_0$ is $3$-stable by
Lemma~\ref{lem:stable-closure}(a). Applying Lemma~\ref{lem:stable-closure}(b) $j$ times, each $4^jN_0$ is $3$-stable,
and so is their union. This proves the case $\kappa=0$.

Let $\kappa\ge1$. By Proposition~\ref{prop:two-recursion-sets},
$\diagSet{2}[\kappa]=B_\kappa\cup4\diagSet{2}[\kappa-1]$, where $B_\kappa$ is either $1+4\padicInts{2}$ or
$1+8\padicInts{2}$. In both cases $B_\kappa$ is a coset $1+2^s\padicInts{2}$ with $s\le3=\val{2}(1)+3$, hence $3$-stable
by Lemma~\ref{lem:stable-closure}(d), and $4\diagSet{2}[\kappa-1]$ is $3$-stable by Lemma~\ref{lem:stable-closure}(b)
and the induction hypothesis. By Lemma~\ref{lem:stable-closure}(a) their union is $3$-stable. This completes the proof.
\end{proof}

\begin{lemma}[Regularity at $p=2$]\label{lem:regularity-two}
For every finite $\shifts\subset\ZZ$ with $\#\shifts=k\ge1$ and every $\alpha\ge1$:
\[
0\ \le\ \truncDen{2}{\alpha}-\locDen{2}[\shifts]\ \le\ \frac{4k}{2^\alpha}.
\]
\end{lemma}

\begin{proof}
The proof is identical to that of Lemma~\ref{lem:regularity-odd}, with Lemma~\ref{lem:stability-two} in place of
Lemma~\ref{lem:stability-odd}. The set $\locSet{2}$ is closed and $3$-stable, so
Lemma~\ref{lem:truncation} with $p=2$ and $c=3$ bounds the single-index discrepancy
$\truncDen{2}{\alpha}[\{0\}]-\locDen{2}$ by $2^{c-1}/2^\alpha=4/2^\alpha$. The reduction to the single-index
discrepancy by the union bound over the $k$ shifts (Claims 1 and 2 in the proof of Lemma~\ref{lem:regularity-odd})
does not use the parity of $p$, and gives the stated bound.
\end{proof}

\subsection{Summary}\label{subsec:local-summary}

Table~\ref{tab:local-summary} collects the results of \S\ref{sec:local-odd} and \S\ref{sec:local-two}. In it,
$\locSet{p}$ is described up to multiplication by a unit, which changes neither the measure nor the stability
(Lemma~\ref{lem:stable-closure}(c)).

\begin{table}[!htbp]
\centering
\begin{tabular}{llll}
\hline
Prime & $\locSet{p}$ (up to a unit) & $\haar{p}(\locSet{p})$ & $c$\\
\hline
$p$ odd, $\kron{p}=1$ & $\padicInts{p}$ & $1$ & $1$\\
$p$ odd, $\kron{p}=-1$ & $\{0\}\cup\{\val{p}\text{ even}\}$ & $p/(p+1)$ & $1$\\
$p$ odd, $p\mid d$ & recursion in $e$ (\S\ref{subsec:ramified}) &
Props.~\ref{prop:ramified-odd}--\ref{prop:ramified-even-inert} & $1$\\
$p=2$, $d\equiv1\bmod8$ & $\padicInts{2}$ & $1$ & $1$\\
$p=2$, $d\equiv5\bmod8$ & $\{0\}\cup\{\val{2}\text{ even}\}$ & $2/3$ & $1$\\
$p=2$, $4\mid d$ & recursion in $\kappa$ (Prop.~\ref{prop:two-recursion-sets}) &
Props.~\ref{prop:two-base-1}--\ref{prop:two-recursion} & $3$\\
\hline
\end{tabular}
\caption{The value sets $\locSet{p}$, their measures, and the constant $c$ for which they are $c$-stable
(Lemmas~\ref{lem:stability-odd} and~\ref{lem:stability-two}).}
\label{tab:local-summary}
\end{table}
\FloatBarrier

\section{The singular series and its average}\label{sec:singular-series}

In this section we show that the product $\singSer=\prod_p\locFac{p}$ converges, and that
$\singSer$ has the correct average over a window of shifts. For the latter we adapt to our setting the
simplified argument of Ford~\cite{ford}, which replaces Gallagher's original proof~\cite{gallagher} for the primes; the
exact local densities and regularity lemmas of \S\ref{sec:local-odd} and \S\ref{sec:local-two} replace the
corresponding facts for the primes.

\subsection{A Fubini/CRT averaging identity}\label{subsec:averaging}

\begin{lemma}[Averaging identity]\label{lem:averaging}
Let $G$ be a finite abelian group, $T\subseteq G$ any subset, $\mu(T)=\#T/\#G$. For $k\ge1$,
\[
\frac1{\#G^k}\sum_{\shifts\in G^k}\mu\Big(\bigcap_i\left(T-h_i\right)\Big)=\mu(T)^k.
\]
\end{lemma}

\begin{proof}
We expand the left-hand side:
\[
\frac1{\#G^k}\sum_{\shifts}\frac1{\#G}\#\{a\in G:\ a+h_i\in T\ \forall i\}
=\frac1{\#G^{k+1}}\sum_{\shifts}\sum_a\prod_i\mathbf 1_T(a+h_i).
\]
Now swap the order of summation. For fixed $a$, we have:
\[
\sum_{\shifts\in G^k}\prod_i\mathbf1_T(a+h_i)=\prod_i\Big(\sum_{h_i\in G}\mathbf1_T(a+h_i)\Big)=(\#T)^k,
\]
since each factor sums to $\#T$ independently, as $h_i\mapsto a+h_i$ is a bijection of $G$. So the total is:
\[
\frac1{\#G^{k+1}}\sum_a(\#T)^k=\frac{\#G\cdot(\#T)^k}{\#G^{k+1}}=\mu(T)^k,
\]
which is exactly what we wanted to show.
\end{proof}

Apply this with $G=\ZZ/p^\alpha\ZZ$ and $T$ the set defining $\truncDen{p}{\alpha}$ at each prime
(Definition~\ref{def:truncated}). Since Lemma~\ref{lem:averaging} holds for \emph{any} subset $T$, it applies verbatim,
and gives:
\[
\frac1{p^{k\alpha}}\sum_{\shifts\bmod p^\alpha}\truncFac{p}{\alpha}=1,
\]
where $\truncFac{p}{\alpha}$ is the truncated local factor of Lemma~\ref{lem:monotonicity}, and $k$ is the
length of the tuple $\shifts$, even when its entries repeat. Thus the truncated local factor averages
to exactly $1$, at every prime and every level, regardless of what $\locSet{p}$ actually is.

\subsection{Convergence of the singular series}\label{subsec:convergence}

\begin{corollary}[Convergence of the singular series]\label{cor:convergence}
For every finite $\shifts\subset\ZZ$ with $\#\shifts=k$, the product
$\singSer=\prod_p\locFac{p}$ converges absolutely, and:
\[
\singSer\ \ll_{k,d}\ \big(\log\log(|\det\shifts|+3)\big)^{k-1}.
\]
If $\locFac{p}=0$ for some $p$, then we set $\singSer=0$, and absolute convergence refers to
the remaining factors.
\end{corollary}

\begin{proof}
At an odd split prime we have $\locSet{p}=\padicInts{p}$ (Proposition~\ref{prop:split}), so $\locFac{p}=1$ and these
primes may be ignored. We split the rest of the product into four kinds of primes, which are clearly disjoint (by
parity, and by the trichotomy of \S\ref{subsec:trichotomy}): $p=2$, the odd primes dividing $d$, the odd inert primes
dividing $\det(\shifts)$ (at most $\omega(\det\shifts)$ of them), and the remaining odd inert primes. The first two
kinds together are exactly the primes dividing $2d$, so there are $r:=\omega(2d)$ of them. We treat three cases.

\begin{enumerate}
\item[Case 1:] The first two kinds, that is, $p=2$ and the odd primes dividing $d$. Here $\locDen{p}\ge1/6$
uniformly. Indeed, Propositions~\ref{prop:ramified-odd}--\ref{prop:ramified-even-inert} give
$\locDen{p}>p/(2(p+1))\ge3/8$ for odd $p\mid d$, while for $p=2$
Propositions~\ref{prop:two-split},~\ref{prop:two-inert}, and~\ref{prop:two-base-1}--\ref{prop:two-base-even} give
$\locDen{2}\ge1/2$ when $2\nmid d$ or $4\mid d$ with conductor depth $\kappa=0$, and
Proposition~\ref{prop:two-recursion} gives $\locDen{2}>1/6$ for every conductor depth $\kappa\ge1$.
Lemma~\ref{lem:monotonicity} therefore gives $\locFac{p}\le\locDen{p}^{1-k}\le6^{k-1}$.

\item[Case 2:] The odd inert primes dividing $\det(\shifts)$. Here Lemma~\ref{lem:universal-near-one} gives
$\locFac{p}\le(1+1/p)^{k-1}$. The number of such primes is at most $\omega(n)$, where $n:=|\det\shifts|$, and is not
bounded in terms of $k$ alone, so we must bound the product of these factors with some care. A product $\prod_{p\in
P}(1+1/p)$ over a set $P$ of $\omega(n)$ primes is largest when $P$ consists of the first $\omega(n)$ primes. Since
$\omega(n)\ll\log n/\log\log n$, the $\omega(n)$-th prime is $\ll\log n$, and Mertens' second theorem gives:
\[
\prod_{\substack{p\mid\det\shifts\\p\text{ odd inert}}}\locFac{p}
\ \le\ \prod_{p\ll\log n}\Big(1+\frac1p\Big)^{k-1}
\ \ll_k\ \big(\log\log(n+3)\big)^{k-1}.
\]

\item[Case 3:] The remaining odd inert primes. Here Lemma~\ref{lem:near-one} gives
$\locFac{p}=1+O_k(p^{-2})$, and $\sum_pp^{-2}<\infty$ closes the argument.
\end{enumerate}

Putting these together, the first two kinds contribute at most $6^{(k-1)r}$, a constant depending only on $k$ and $d$;
the third kind contributes $\ll_k(\log\log(|\det\shifts|+3))^{k-1}$, a finite product for each fixed $\shifts$; and
the last kind converges absolutely to a quantity $\ll_k1$. This completes the proof.
\end{proof}

\subsection{The Ford-style average}\label{subsec:ford}

\begin{theorem}[Average of the singular series]\label{thm:ford-average}
Fix $k\ge1$. As $H\to\infty$, with the sum running over integer tuples:
\[
\sum_{0<h_1<\cdots<h_k<H}\singSer=\frac{H^k}{k!}\Big(1+O_{k,d}\Big(\frac1{\log\log H}\Big)\Big).
\]
\end{theorem}

The core idea of the proof, following Ford~\cite{ford}, is to truncate the Euler product at a cutoff $Y=c\log H$, where
$0<c<1/2$ is fixed. The primes above $Y$ contribute a factor $1+o(1)$, uniformly in $\shifts$, while the primes up to
$Y$ are replaced by a periodic function $\Phi$ whose period $Q$ is at most $H^{2c+o(1)}$, a power of $H$ strictly
smaller than $H$, and whose sum over a full period is computed exactly from Lemma~\ref{lem:averaging} via the Chinese
remainder theorem.

\begin{proof}
Recall from the proof of Corollary~\ref{cor:convergence} that $p=2$ together with the odd primes dividing $d$, that is,
the $r=\omega(2d)$ primes dividing $2d$, satisfy $\locDen{p}\ge1/6$ uniformly. We call these the $r$
\emph{exceptional} primes below. Fix $0<c<1/2$, and set $Y:=c\log H$; for $H$ large enough (depending only on $d$), $Y$
exceeds all $r$ exceptional primes, so that only odd inert primes contribute to $p>Y$. Set $L:=\log H\,(\log\log H)^k$,
and for each $p\le Y$ set $\alpha_p:=\lceil\log L/\log p\rceil$, so that:
\[
L\le p^{\alpha_p}<pL.
\]
Let $Q:=\prod_{p\le Y}p^{\alpha_p}$, and define the $Q$-periodic function $\Phi(\shifts):=\prod_{p\le
Y}\truncFac{p}{\alpha_p}$. Throughout, implied constants may depend on $k$ and $d$. Let:
\[
D:=\{\shifts\in[0,H]^k\cap\ZZ^k:\ h_1,\dots,h_k\text{ pairwise distinct}\}.
\]
The singular series $\singSer$ is only defined for $\shifts\in D$ (Definition~\ref{def:singular}), and
indeed $\det\shifts=0$ for every tuple with a repeated coordinate, whereas the finite product $\Phi(\shifts)$ is
defined for every $\shifts\in\ZZ^k$. We split the proof into five claims; Claims 1 and 4 concern only tuples in $D$,
while Claim 5 is a statement about the full box.

Claim 1: We have, uniformly in $\shifts\in D$:
\[
\prod_{\substack{p>Y\\\kron{p}=-1}}\locFac{p}=1+O_k\Big(\frac1{\log\log H}\Big).
\]

Proof of Claim 1: We split by whether $p\mid\det(\shifts)$. For $p\nmid\det(\shifts)$, Lemma~\ref{lem:near-one}
gives $\locFac{p}=1+O_k(p^{-2})$, and $\sum_{p>Y}p^{-2}=O(1/Y)$. For $p\mid\det(\shifts)$,
Lemma~\ref{lem:universal-near-one} gives the upper bound $\locFac{p}\le(1+1/p)^{k-1}=1+O_k(1/p)$. For the
lower bound, observe that every unit lies in $\locSet{p}$ (Proposition~\ref{prop:inert}), so $\locDen{p}[\shifts]$ is at
least the measure of the set of $a\in\padicInts{p}$ for which every $a+h_i$ is a unit. This set is the complement of at
most $k$ residue classes modulo $p$, so $\locDen{p}[\shifts]\ge1-k/p$. Since $\locDen{p}\le1$, we get
$\locFac{p}\ge1-k/p$. Together these give $\locFac{p}=1+O_k(1/p)=1+O_k(1/Y)$, and there are at
most $\omega(\det\shifts)=O_k(\log H/\log\log H)$ such primes, since $\shifts\in D$ gives
$0<|\det\shifts|\le H^{\binom k2}$ and the number of distinct prime factors of an integer $n\ge1$ is
$O(\log n/\log\log n)$. So we obtain:
\begin{align*}
\prod_{\substack{p>Y\\\kron{p}=-1}}\locFac{p}
&=\big(1+O_k(1/Y)\big)\cdot\big(1+O_k(1/Y)\big)^{O_k(\log H/\log\log H)}\\
&=1+O_k\Big(\frac{\log H}{Y\log\log H}\Big),
\end{align*}
uniformly in $\shifts\in D$. Substituting $Y=c\log H$, the last expression is $1+O_k(1/\log\log H)$. This proves
Claim 1.

Claim 2: We have $Q\le H^{2c+o(1)}$. In particular, setting $\eta:=(1-2c)/2>0$, we have $Q\le H^{1-\eta}$ for all
sufficiently large $H$.

Proof of Claim 2: Since $p^{\alpha_p}<pL$ for every $p\le Y$, we have:
\[
\log Q=\sum_{p\le Y}\alpha_p\log p<\sum_{p\le Y}(\log L+\log p)=\pi(Y)\log L+\vartheta(Y),
\]
where $\vartheta(Y):=\sum_{p\le Y}\log p$. By the prime number theorem, $\pi(Y)\sim Y/\log Y$ and $\vartheta(Y)\sim Y$.
Moreover $\log L=\log\log H+k\log\log\log H\sim\log\log H\sim\log Y$, so $\pi(Y)\log L\sim Y$. Putting these together
we get:
\[
\log Q<(2+o(1))\,Y=(2c+o(1))\log H,
\]
that is, $Q\le H^{2c+o(1)}$. Since $2c<1-\eta$, the second statement follows. This proves Claim 2.

Claim 3: We have $\|\Phi\|_\infty\ll_{k,d}(\log\log H)^{k-1}$, and the identical bound holds for
$\prod_{p\le Y}\locFac{p}$, uniformly in $\shifts\in D$.

Proof of Claim 3: At an odd split prime we have $\locSet{p}=\padicInts{p}$ (Proposition~\ref{prop:split}), so every
factor at $p$, whether truncated or exact, equals $1$. We split the remaining primes $p\le Y$ into the $r$ exceptional
primes and the odd inert primes $\le Y$. For the exceptional primes, the lower bound in Lemma~\ref{lem:regularity-odd}
(respectively Lemma~\ref{lem:regularity-two} at $p=2$) gives $\truncDen{p}{\alpha_p}[\{0\}]\ge\locDen{p}\ge1/6$, so
$\truncFac{p}{\alpha_p}\le\truncDen{p}{\alpha_p}[\{0\}]^{1-k}\le6^{k-1}$ by Lemma~\ref{lem:monotonicity},
and these $r$ primes contribute at most $6^{(k-1)r}$ to $\Phi$, a constant depending only on $k,d$.

For the odd inert primes $\le Y$, Proposition~\ref{prop:inert} gives $\locDen{p}=p/(p+1)$ exactly, and the same
lower bound in Lemma~\ref{lem:regularity-odd} gives $\truncDen{p}{\alpha_p}[\{0\}]\ge p/(p+1)$. So
Lemma~\ref{lem:monotonicity} gives $\truncFac{p}{\alpha_p}\le(1+1/p)^{k-1}$, which is identical to the
untruncated bound of Lemma~\ref{lem:universal-near-one}. We find:
\begin{align*}
\prod_{\substack{p\le Y\\p\text{ odd inert}}}\truncFac{p}{\alpha_p}
&\ \le\ \prod_{p\le Y}\Big(1+\frac1p\Big)^{k-1}\ \ll_k\ \exp\Big((k-1)\sum_{p\le Y}\frac1p\Big)\\
&\ \stackrel{(\ast)}{\ll_k}\ (\log Y)^{k-1}\ \ll_k\ (\log\log H)^{k-1},
\end{align*}
where $(\ast)$ is Mertens' second theorem, and the last step uses $\log Y\sim\log\log H$. Combining the two kinds of
primes we get $\|\Phi\|_\infty\ll_{k,d}(\log\log H)^{k-1}$, and the identical bound holds for
$\prod_{p\le Y}\locFac{p}$ by the same computation with $\alpha_p$ dropped. This proves Claim 3.

Claim 4: We have $\Big|\prod_{p\le Y}\locFac{p}-\Phi(\shifts)\Big|\ll(\log\log H)^{-2}$, uniformly in
$\shifts\in D$.

Proof of Claim 4: For $p\le Y$, write $a:=\locDen{p}[\shifts]$, $a':=\truncDen{p}{\alpha_p}$,
$b:=\locDen{p}$, and $b':=\truncDen{p}{\alpha_p}[\{0\}]$. By Lemma~\ref{lem:regularity-odd} (or
Lemma~\ref{lem:regularity-two} at $p=2$), we have $0\le a'-a\le4k/p^{\alpha_p}$ and $0\le b'-b\le4/p^{\alpha_p}$ (the
constant $4$ is needed only at $p=2$), and $b,b'\ge1/6$ throughout, by the proof of Claim 3. Moreover $a'\le b'\le1$
(Lemma~\ref{lem:monotonicity}), and $|1/b^k-1/b'^k|\le k\max(1/b,1/b')^{k+1}(b'-b)$ by the mean value bound. This
gives:
\begin{align*}
\big|\locFac{p}-\truncFac{p}{\alpha_p}\big|
&=\Big|\frac a{b^k}-\frac{a'}{b'^k}\Big|\\
&\stackrel{(\ast)}{\le}\ \frac{a'-a}{b^k}+a'\Big|\frac1{b^k}-\frac1{b'^k}\Big|\\
&\stackrel{(\ast\ast)}{\le}\ \frac{4k\cdot6^k}{p^{\alpha_p}}+\frac{4k\cdot6^{k+1}}{p^{\alpha_p}}\ =\
\frac{C_k}{p^{\alpha_p}},
\end{align*}
for a constant $C_k$ depending only on $k$. Here $(\ast)$ follows from the triangle inequality after adding and
subtracting $a'/b^k$, and $(\ast\ast)$ from the bounds on $a'-a$, $b,b'$, and $a'$ just recorded, together with the mean
value bound.

To pass from this per-prime bound to a bound on the difference of the two full products, write $p_1<\cdots<p_M$ for the
primes $\le Y$ ($M=\pi(Y)$), and for $0\le j\le M$ define the \emph{hybrid} product with the first $j$ factors truncated
and the rest exact:
\[
\Pi_j:=\truncFac{p_1}{\alpha_{p_1}}\cdots\truncFac{p_j}{\alpha_{p_j}}
\cdot\locFac{p_{j+1}}\cdots\locFac{p_M}.
\]
Then $\Pi_0=\prod_{p\le Y}\locFac{p}$ and $\Pi_M=\Phi(\shifts)$, so telescoping gives:
\[
\prod_{p\le Y}\locFac{p}-\Phi(\shifts)\ =\ \Pi_0-\Pi_M\ =\ \sum_{j=1}^M(\Pi_{j-1}-\Pi_j).
\]
Since $\Pi_{j-1},\Pi_j$ share every factor except the $j$-th (exact in $\Pi_{j-1}$, truncated in $\Pi_j$), we have:
\[
\Pi_{j-1}-\Pi_j\ =\ \Big(\prod_{i<j}\truncFac{p_i}{\alpha_{p_i}}\cdot
\prod_{i>j}\locFac{p_i}\Big)
\cdot\big(\locFac{p_j}-\truncFac{p_j}{\alpha_{p_j}}\big).
\]
To bound the product of the $M-1$ surrounding factors, set $c_p:=6^{k-1}$ if $p$ is one of the $r$ exceptional primes,
$c_p:=(1+1/p)^{k-1}$ if $p$ is odd inert, and $c_p:=1$ if $p$ is odd split. By the proof of Claim 3, every factor at
$p$, whether truncated or exact, satisfies $0\le\locFac{p}\le c_p$ and $0\le\truncFac{p}{\alpha_p}\le c_p$, the lower
bound holding because each factor is a ratio of densities. Hence the product of the surrounding factors is at most
$\prod_{i\ne j}c_{p_i}$, and since $c_{p_j}\ge1$ we get:
\[
\prod_{i<j}\truncFac{p_i}{\alpha_{p_i}}\cdot\prod_{i>j}\locFac{p_i}
\ \le\ \prod_{i\ne j}c_{p_i}\ \le\ \prod_{p\le Y}c_p\ \ll_{k,d}\ (\log\log H)^{k-1},
\]
where the last bound is the computation in the proof of Claim 3. Combined with the per-prime bound above,
$|\Pi_{j-1}-\Pi_j|\ll_{k,d}(\log\log H)^{k-1}\cdot C_k/p_j^{\alpha_{p_j}}$, and summing over $j=1,\dots,M$ we find:
\begin{align*}
\Big|\prod_{p\le Y}\locFac{p}-\Phi(\shifts)\Big|
&\le\ \sum_{j=1}^M|\Pi_{j-1}-\Pi_j|\ \ll_{k,d}\ (\log\log H)^{k-1}\sum_{p\le Y}p^{-\alpha_p}\\
&\le\ (\log\log H)^{k-1}\cdot\frac{\pi(Y)}L\ \ll\ \frac{(\log\log H)^{k-1}}{(\log\log H)^{k+1}}
\ =\ \frac1{(\log\log H)^2},
\end{align*}
uniformly in $\shifts \in D$, using $p^{\alpha_p}\ge L$ and $\pi(Y)\ll\log H/\log\log H$ (by the prime number theorem,
as in the proof of Claim 2). This proves Claim 4.

Claim 5: We have $\sum_{\shifts\in[0,H]^k\cap\ZZ^k}\Phi(\shifts)=H^k+O\big(H^{k-\eta}(\log\log H)^{k-1}\big)$.

Proof of Claim 5: Since $\Phi$ is $Q$-periodic and $Q\le H^{1-\eta}$ by Claim 2, we may group $[0,H]^k\cap\ZZ^k$ by
residue class mod $Q$: each class contains $(H/Q)^k+O((H/Q)^{k-1})$ lattice points. By Lemma~\ref{lem:averaging},
applied at each $p\le Y$ separately and combined via the Chinese remainder theorem (the moduli $p^{\alpha_p}$, $p\le Y$,
being pairwise coprime), we get:
\[
\sum_{\mathbf a\bmod Q}\Phi(\mathbf a)=\prod_{p\le Y}\sum_{\mathbf a\bmod p^{\alpha_p}}
\truncFac{p}{\alpha_p}[\mathbf a]
=\prod_{p\le Y}p^{k\alpha_p}=Q^k.
\]
Hence, using Claims 2 and 3:
\begin{align*}
\sum_{\shifts\in[0,H]^k\cap\ZZ^k}\Phi(\shifts)
&=\Big(\frac HQ\Big)^kQ^k+O\Big(\Big(\frac HQ\Big)^{k-1}\|\Phi\|_\infty\,Q^k\Big)\\
&=H^k+O\big(H^{k-1}Q(\log\log H)^{k-1}\big)=H^k+O\big(H^{k-\eta}(\log\log H)^{k-1}\big),
\end{align*}
where we used $Q\le H^{1-\eta}$. This proves Claim 5.

We now put all of our claims to work. The tuples in $[0,H]^k\cap\ZZ^k$ with a repeated coordinate number
$O(H^{k-1})$, and on each of them $0\le\Phi(\shifts)\le\|\Phi\|_\infty$, so by Claims 5 and 3:
\[
\sum_{\shifts\in D}\Phi(\shifts)=H^k+O\big(H^{k-\eta}(\log\log H)^{k-1}\big)+O\big(H^{k-1}(\log\log H)^{k-1}\big)
=H^k\big(1+O(H^{-\eta/2})\big).
\]
From now on every tuple lies in $D$, where $\singSer$ is defined and Claims 1 and 4 apply. By Claim 4,
replacing $\Phi(\shifts)$ by $\prod_{p\le Y}\locFac{p}$ changes each of the $O(H^k)$ terms by
$O((\log\log H)^{-2})$, so that:
\[
\sum_{\shifts\in D}\prod_{p\le Y}\locFac{p}=H^k\Big(1+O\Big(\frac1{(\log\log H)^2}\Big)\Big).
\]
Every term of this sum is nonnegative, and multiplying it by the tail factor of Claim 1, which is
$1+O(1/\log\log H)$ uniformly in $\shifts\in D$, gives $\singSer$. Hence:
\[
\sum_{\shifts\in D}\singSer=H^k\Big(1+O\Big(\frac1{\log\log H}\Big)\Big).
\]
Finally, each $k$-subset of $(0,H)\cap\ZZ$ is counted $k!$ times among the tuples in $D$, once for each ordering. The
remaining tuples in $D$, those with a coordinate equal to $0$ or to $H$, number $O(H^{k-1})$, and by Claims 1 and 3 each
of them has $\singSer\ll(\log\log H)^{k-1}$. We conclude:
\[
\sum_{0<h_1<\cdots<h_k<H}\singSer=\frac{H^k}{k!}\Big(1+O\Big(\frac1{\log\log H}\Big)\Big),
\]
where the implied constant depends on $k$ and $d$. This completes the proof.
\end{proof}

\section{\texorpdfstring{Proof of Theorem~\ref{thm:genus}}{Proof of the genus-level theorem}}\label{sec:proof-genus}

In this section we prove Theorem~\ref{thm:genus} by the method of moments, loosely following Gallagher~\cite{gallagher}:
under Hypothesis~$U(k)$, every factorial moment of the number of elements of $E$ in a short window converges to the
corresponding moment of a Poisson distribution.

Fix $\tau>0$. By \S\ref{sec:bernays} (Inputs~\ref{inp:asymptotic} and~\ref{inp:equidistribution} combined), we have
$\repDensity{x}=N(x)/x\sim \bConst{d}/\sqrt{\log x}$. By the hypothesis of Theorem~\ref{thm:genus}, the window
$y=y(x,\tau)$ satisfies:
\[
y\cdot \repDensity{x}\longrightarrow\tau\qquad(\text{so } y\sim\tau\sqrt{\log x}/\bConst{d}),
\]
and, for $n$ uniform on the integers $1\le n\le x$, we let $Y_x$ count $E\cap(n,n+y]$, that is, $Y_x=\sum_{0<t\le
y}\mathbf1_E(n+t)$. Note that $y\asymp\sqrt{\log x}$, so that $(0,y]\subset(0,C\sqrt{\log x}\,]$ for a suitable constant
$C$, which is a range covered by Hypothesis~$U(k)$.

\begin{theorem}[Factorial moments, conditional on Hypothesis $U(k)$]\label{thm:moments}
For every fixed $k\ge1$, the $k$-th factorial moment $\mathbb E_n[(Y_x)_k]\to\tau^k$ as $x\to\infty$.
\end{theorem}

The idea of the proof is to write the factorial moment as a sum of $k$-point correlations over $k$-tuples of shifts in
the window, to replace each $k$-point correlation by the singular series using Hypothesis~$U(k)$, and then to evaluate
the resulting sum with Theorem~\ref{thm:ford-average}.

\begin{proof}
We start by expanding the factorial moment as a sum over ordered $k$-tuples of distinct shifts in $(0,y]$. Since $n$ is
uniform on the $\lfloor x\rfloor$ integers $1\le n\le x$, while $\kCorr{k}{x}$ is normalized by $x$, we get:
\[
\mathbb E_n[(Y_x)_k]=\frac x{\lfloor x\rfloor}\,M_k(x),\qquad\text{where}\quad
M_k(x):=k!\sum_{0<h_1<\cdots<h_k\le y}\kCorr{k}{x}.
\]
Since $x/\lfloor x\rfloor\to1$, it suffices to show that $M_k(x)\to\tau^k$. By Hypothesis~$U(k)$, uniformly over this
range, $\kCorr{k}{x}=\singSer\repDensity{x}^k+o(\repDensity{x}^k)$. The error term is additive, since
$\singSer$ may vanish. There are at most $y^k$ tuples in the range, so summing gives:
\[
M_k(x)=k!\,\repDensity{x}^k\sum_{0<h_1<\cdots<h_k\le y}\singSer\ +\ o\big((y\repDensity{x})^k\big).
\]
This range agrees with that of Theorem~\ref{thm:ford-average} with $H=y$, except possibly for tuples with $h_k=y$.
There are $O(y^{k-1})$ such tuples, each with $|\det\shifts|\le y^{\binom k2}$ and hence contributing
$\singSer\ll_{k,d}(\log\log y)^{k-1}$ by Corollary~\ref{cor:convergence}, so that they contribute
$O_{k,d}(y^{k-1}(\log\log y)^{k-1})=o(y^k)$ in total. By Theorem~\ref{thm:ford-average} with $H=y$, the sum is therefore
$y^k(1/k!+o(1))$, so that:
\[
M_k(x)=k!\cdot \repDensity{x}^k\cdot y^k\Big(\frac1{k!}+o(1)\Big)+o\big((y\repDensity{x})^k\big)
=\big(y\repDensity{x}\big)^k(1+o(1))\longrightarrow\tau^k,
\]
directly from $y\repDensity{x}\to\tau$. This completes the proof.
\end{proof}

\begin{proof}[Proof (of Theorem~\ref{thm:genus})]
The factorial moments of $\mathrm{Poisson}(\tau)$ are exactly $\tau^k$, so by Theorem~\ref{thm:moments} every factorial
moment of $Y_x$, and hence every ordinary moment, converges to that of $\mathrm{Poisson}(\tau)$. Since the Poisson
distribution is determined by its moments, the method of moments~\cite[\S30]{billingsley} gives
$Y_x\to\mathrm{Poisson}(\tau)$ in distribution. Since $Y_x=N(n+y)-N(n)$ takes integer values, this means precisely that
for every fixed $m\ge0$:
\[
\mathbb P(Y_x=m)\ =\ \frac1{\lfloor x\rfloor}\,\#\big\{1\le n\le x:\ N(n+y)-N(n)=m\big\}\ \longrightarrow\
e^{-\tau}\,\frac{\tau^m}{m!}.
\]
Since $\lfloor x\rfloor/x\to1$, the same limit holds with $\lfloor x\rfloor$ replaced by $x$, which is the statement of
Theorem~\ref{thm:genus}. This completes the proof.
\end{proof}

Combined with \S\ref{sec:class-from-genus}, this also completes the proof of Theorem~\ref{thm:class}.

\section{\texorpdfstring{Numerical evidence for Hypothesis $U(k)$}{Numerical evidence for Hypothesis U(k)}}
\label{sec:numerics}

In this section we test Hypothesis~$U(k)$ numerically for $k=2,3$, by comparing the normalized $k$-point correlations
$\kCorr{k}{x}/\repDensity{x}^k$ of $\repSet{R}$ (Definition~\ref{def:k-point-correlation}) with the singular series
$\singSer$ computed from the exact local densities of \S\ref{sec:local-odd} and \S\ref{sec:local-two},
and then by looking directly at the window counts. The computations were carried out in the Wolfram Language.

\subsection{Method}\label{subsec:numerics-method}

We took $d\in\{-3,-4,-7,-8,-12,-16,-20\}$, with both genera for $d=-20$ (those of $x^2+5y^2$ and $2x^2+2xy+3y^2$); the
discriminants $-12$ and $-16$ are not fundamental, and $d=-4$ is the case of sums of two squares studied
in~\cite{freiberg-kurlberg-rosenzweig}. For each genus $R$ we sieved the indicator function of $E=\repSet{R}$ up to
$x=10^7$ from the values of every reduced form in $R$, and computed $\kCorr{k}{x}/\repDensity{x}^k$ exactly, taking
$h_1=0$. This is harmless: $\singSer$ is invariant under translating $\shifts$, while translating $\shifts$ by
$-h_1$ only shifts the range of $n$ in Definition~\ref{def:k-point-correlation} by $h_1$, so it changes
$\kCorr{k}{x}$ by at most $h_1/x$, and hence the quantity in Hypothesis~$U(k)$ by
$O\big(h_1/(x\repDensity{x}^k)\big)=O\big(h_1(\log x)^{O(1)}/x\big)$.

We evaluated $\singSer=\prod_p\locFac{p}$ from the local densities. Odd split primes
contribute $1$, and for odd inert primes $p\nmid\det(\shifts)$ we used the exact values $\locDen{p}[\shifts]=1-k/(p+1)$
and $\locDen{p}=p/(p+1)$ (Lemma~\ref{lem:near-one} and Proposition~\ref{prop:inert}). At the remaining primes, that is
$p=2$, the primes dividing $d$, and the inert primes dividing $\det(\shifts)$, we computed the truncated densities
$\truncDen{p}{\alpha}$ by enumerating residues modulo $p^\alpha\ge2000$, which by Lemmas~\ref{lem:regularity-odd}
and~\ref{lem:regularity-two} is accurate to within $O_k(1/2000)$. Primes up to $10^5$ were included, so the neglected
tail is $O(k^2/(10^5\log10^5))$.

For $k=2$ we took all shifts $(0,h)$ with $1\le h\le60$, and for $k=3$ we sampled $40$ random triples $(0,h_2,h_3)$ with
$0<h_2<h_3\le30$. Recall from \S\ref{sec:k-point-correlations} that a tuple is inadmissible if $\singSer=0$, that is, if
$\locDen{p}[\shifts]=0$ for some prime $p$, as for $\shifts=\{0,1,2\}$ and $d=-3$; no $n$ then has all of $n+h_i$ in
$E$. Such tuples are excluded from the relative errors below. Hypothesis~$U(k)$ concerns shifts up to $C\sqrt{\log x}$,
and $\sqrt{\log x}$ is only about $4$ at $x=10^7$, so our shifts up to $60$ correspond to rather large values of $C$.

\subsection{\texorpdfstring{$2$-point and $3$-point correlations}{2-point and 3-point correlations}}
\label{subsec:numerics-correlations}

\begin{figure}[!htbp]
\centering
\includegraphics[width=0.82\textwidth]{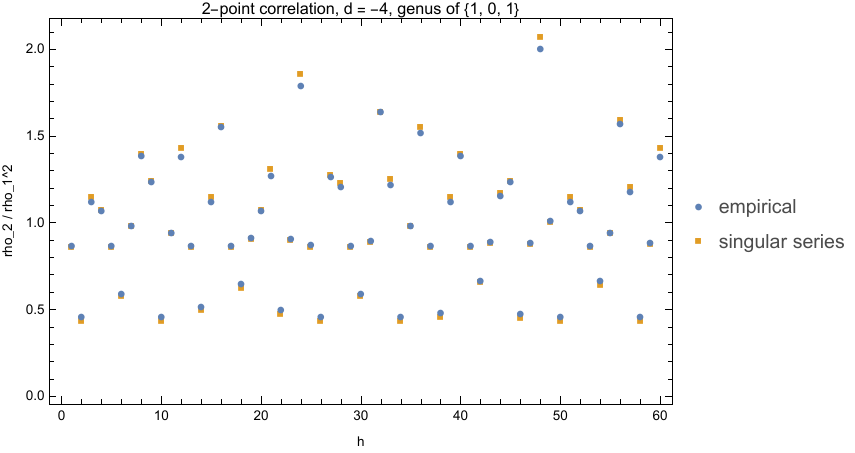}\\[1ex]
\includegraphics[width=0.82\textwidth]{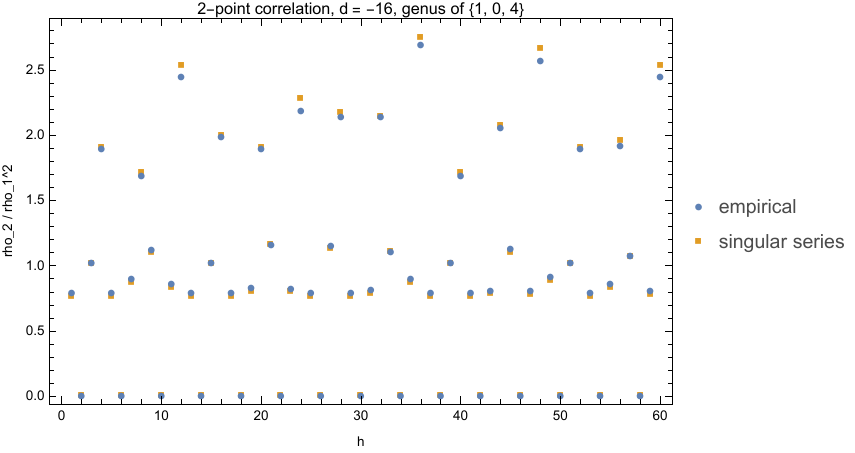}
\caption{$2$-point correlation $\kCorr{2}{x}[\{0,h\}]/\repDensity{x}^2$ at $x=10^7$ (blue dots) against the singular
series $\singSer[\{0,h\}]$ (orange squares), for $d=-4$ (top) and $d=-16$ (bottom).}
\label{fig:pair}
\end{figure}

Figure~\ref{fig:pair} shows that the $2$-point correlations follow the singular series closely through all of its
oscillations. For $d=-3$ (not shown) the largest relative error over $h\le20$ is $5.9\%$, and the mean relative error
over the admissible $h\le60$ is $3.1\%$; note that this range is wider than the one used for $\eps(x)$ in
\S\ref{subsec:numerics-conv}. For $d=-16$ the singular series vanishes at every $h\equiv2\pmod4$, and the empirical
$2$-point correlation is indeed $0$ there. The remaining discrepancies are systematic: the singular series is slightly
too high at its peaks and slightly too low in its troughs, as if the empirical $2$-point correlation were a little
closer to $1$. The other discriminants behave in the same way. The data for $d=-12$ coincide with those for $d=-3$,
because $x^2+3y^2$ and $x^2+xy+y^2$ represent the same integers; this is still a consistency check, since the local
factors at $2$ and $3$ are computed separately for the two discriminants.

The two genera of discriminant $-20$ are instructive. They have the same $2$-point singular series: at $p=2$ and $p=5$
their local value sets are, up to the point $0$, complements of each other, each of measure $1/2$, which forces the
local factors to agree, and at every other prime the two sets coincide. Their empirical $2$-point correlations differ by
a few percent at $x=10^7$ (at $h=1$ they are $0.817$ and $0.757$, against $\singSer[\{0,1\}]\approx0.78$, and at
$h=8$ they are $1.210$ and $1.302$, against $\singSer[\{0,8\}]\approx1.27$), but this looks like a lower-order
effect: the mean absolute difference over $h\le20$ decreases from $0.116$ at $x=10^4$ to $0.052$ at $x=10^7$. The genus
does affect the singular series for $k=3$. The values of $x^2+5y^2$ are $\equiv0,\pm1\pmod5$ and those of
$2x^2+2xy+3y^2$ are $\equiv0,\pm2\pmod5$, so $\{0,1,2\}$ is inadmissible only for the second genus and $\{0,1,3\}$ only
for the first. Up to $10^7$, the tuple $\{0,1,2\}$ occurs $10{,}561$ times for $x^2+5y^2$ and never for $2x^2+2xy+3y^2$,
while $\{0,1,3\}$ occurs $11{,}207$ times for $2x^2+2xy+3y^2$ and never for $x^2+5y^2$.

\begin{figure}[!htbp]
\centering
\includegraphics[width=0.48\textwidth]{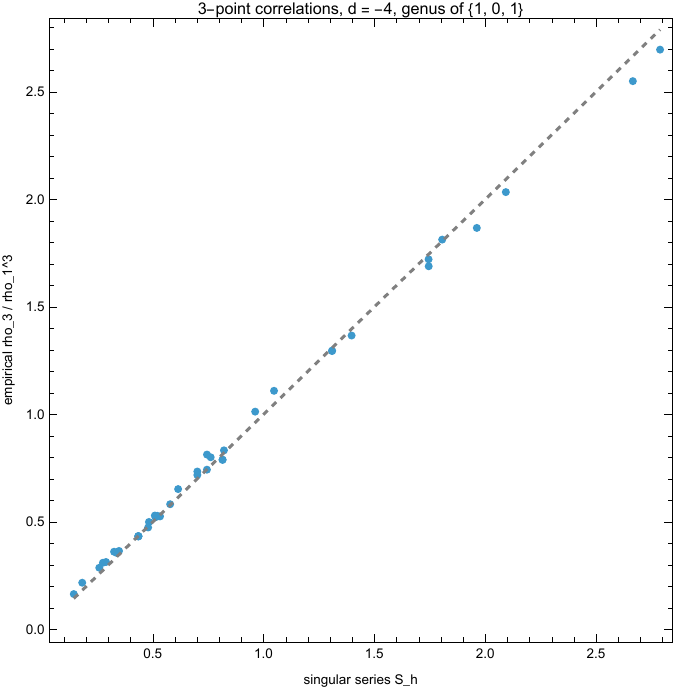}\hfill
\includegraphics[width=0.48\textwidth]{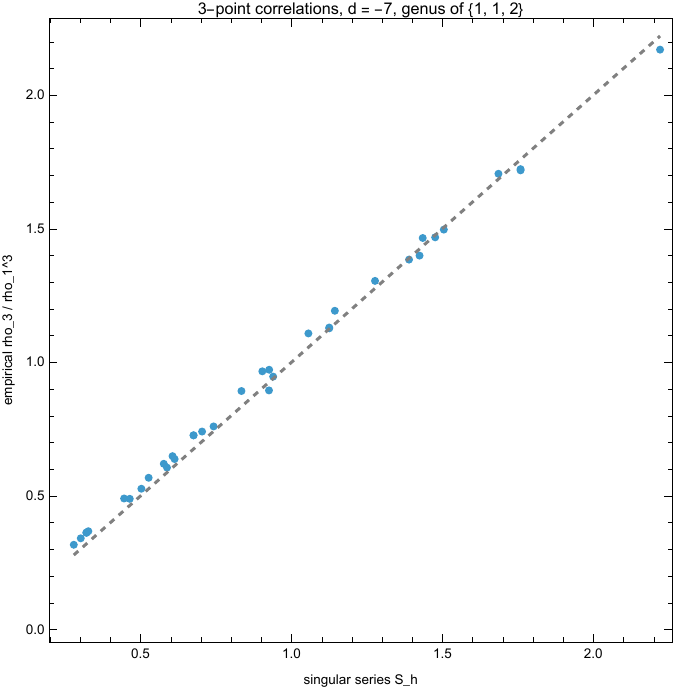}
\caption{$3$-point correlations $\kCorr{3}{x}[\{0,h_2,h_3\}]/\repDensity{x}^3$ at $x=10^7$ against
$\singSer[\{0,h_2,h_3\}]$ for $40$ random triples, for $d=-4$ (left) and $d=-7$ (right). The dashed line is the
identity.}
\label{fig:triple}
\end{figure}

Figure~\ref{fig:triple} shows that the $3$-point correlations lie close to the diagonal over the whole range of values
of $\singSer$, with deviations of a few percent and again slightly below the diagonal at the largest
values. The same holds for the other discriminants.

\subsection{Averages and the dependence on \texorpdfstring{$x$}{x}}\label{subsec:numerics-conv}

\begin{figure}[!htbp]
\centering
\includegraphics[width=0.48\textwidth]{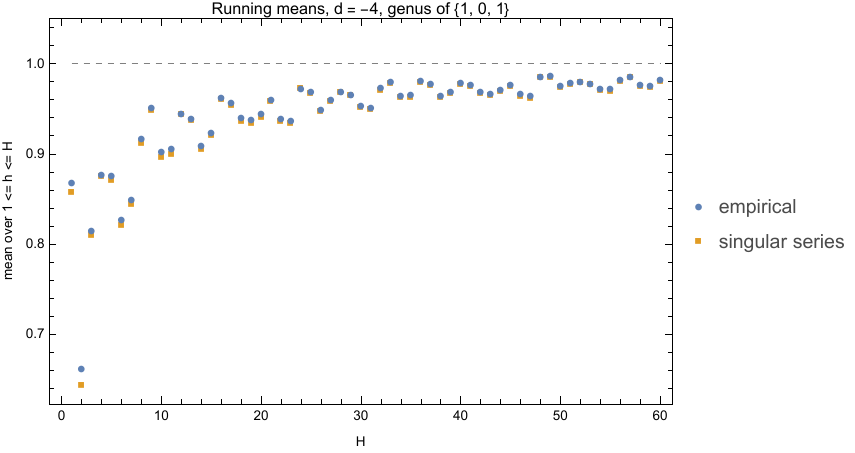}\hfill
\includegraphics[width=0.48\textwidth]{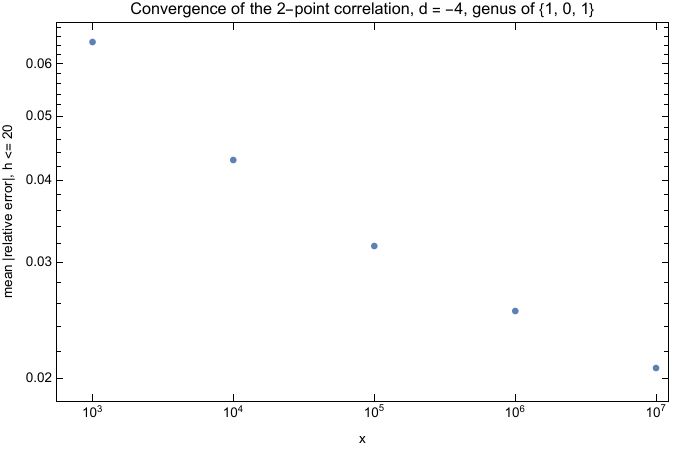}
\caption{Left: the running means $\frac1H\sum_{h=1}^H\kCorr{2}{x}[\{0,h\}]/\repDensity{x}^2$ (blue dots) and
$\frac1H\sum_{h=1}^H\singSer[\{0,h\}]$ (orange squares) for $1\le H\le60$, with $d=-4$ and $x=10^7$; the dashed line
is the constant $1$. Right: the mean relative error $\eps(x)$ of the $2$-point correlation, defined
in~\eqref{eq:mean-rel-error}, for $d=-4$ and $x=10^3,10^4,10^5,10^6,10^7$, on logarithmic axes.}
\label{fig:mean-conv}
\end{figure}

By Theorem~\ref{thm:ford-average} with $k=2$, together with a short monotonicity argument (using $\singSer\ge0$), the
running means $\frac1H\sum_{g\le H}\singSer[\{0,g\}]$ tend to $1$. The left panel of Figure~\ref{fig:mean-conv} shows
that the running means of both the singular series and the empirical $2$-point correlation do tend to $1$, but slowly:
for $d=-4$ they are about $0.88$ at $H=5$, $0.90$ at $H=10$, and $0.98$ at $H=60$, and the two sequences are almost
indistinguishable.

The right panel shows the dependence on $x$. Let $A$ be the set of admissible shifts $1\le h\le20$, that is, those
with $\singSer[\{0,h\}]\ne0$. We measure the agreement at $x$ by the mean relative error:
\begin{equation}\label{eq:mean-rel-error}
\eps(x):=\frac1{|A|}\sum_{h\in A}\bigg|\frac{\kCorr{2}{x}[\{0,h\}]}{\repDensity{x}^2\,\singSer[\{0,h\}]}-1\bigg|.
\end{equation}
For $d=-4$ the mean relative error $\eps(x)$ falls from about $6.5\%$ at $x=10^3$ to about $2\%$ at $x=10^7$; for $d=-3$
it falls from about $9.5\%$ to about $3\%$, and for the genus of $2x^2+2xy+3y^2$ at $d=-20$ from about $12\%$ to about
$4\%$. This is much slower than the $x^{-1/2}$ expected from sampling noise alone, which is consistent with the
logarithmic rate at which $\repDensity{x}=N(x)/x\sim \bConst{d}/\sqrt{\log x}$ settles down: for $d=-4$ the quantity
$N(x)\sqrt{\log x}/x$ is still drifting, from $0.87$ at $x=10^3$ to $0.80$ at $x=10^7$, towards the Landau--Ramanujan
constant $0.7642\ldots$.

\subsection{The distribution in short windows}\label{subsec:numerics-poisson}

\begin{figure}[!htbp]
\centering
\includegraphics[width=0.48\textwidth]{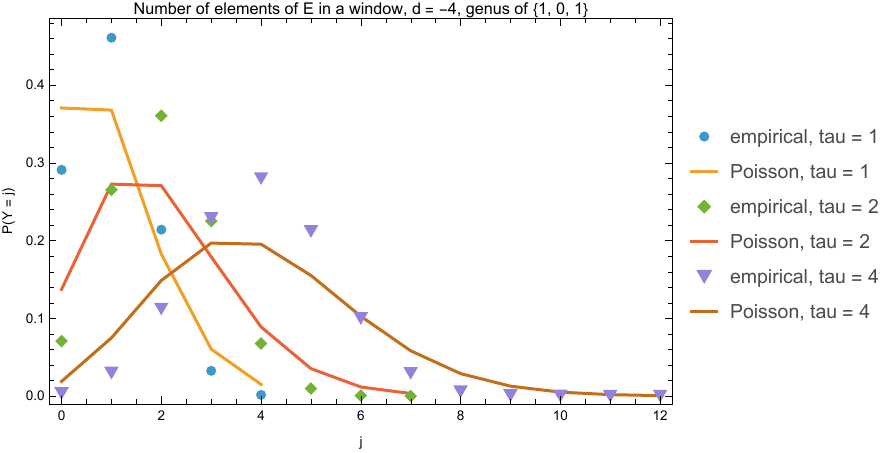}\hfill
\includegraphics[width=0.48\textwidth]{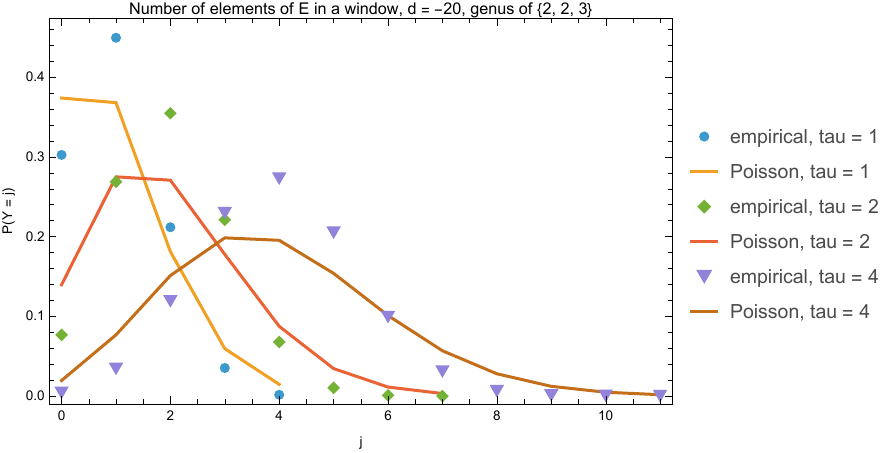}
\caption{Distribution of the number $Y_x$ of elements of $E$ in a window $(n,n+y]$, $n\le x=10^7$, with
$y\repDensity{x}\approx\tau$ for $\tau=1,2,4$ (points), against the Poisson distribution with the same mean (lines), for
$d=-4$ (left) and for the genus of $2x^2+2xy+3y^2$ at $d=-20$ (right).}
\label{fig:poisson}
\end{figure}

Finally, Figure~\ref{fig:poisson} examines the conclusion of Theorem~\ref{thm:genus} directly. At $x=10^7$ the
distribution is clearly not Poisson but under-dispersed: for $\tau=1$ and $d=-4$ the windows contain no element of $E$
in about $29\%$ of cases instead of $e^{-1}\approx37\%$, and exactly one element in about $46\%$ instead of $37\%$; the
genus of $2x^2+2xy+3y^2$ at $d=-20$ gives $30\%$ and $45\%$.

This does not contradict the theorem. The window length $y\sim\tau\sqrt{\log x}/\bConst{d}$ grows only like $\sqrt{\log
x}$, and at $x=10^7$ it is about $5$ for $\tau=1$ and about $20$ for $\tau=4$ when $d=-4$. By the proof of
Theorem~\ref{thm:moments}, the factorial moments of $Y_x$ are governed by the sums of $\singSer$ over the
actual window, and by Figure~\ref{fig:mean-conv} such sums are still noticeably below their limiting value at these
lengths, especially because the small shifts, where $\singSer[\{0,h\}]$ is far from $1$, carry more weight (a window
$(n,n+y]$ contains $\lfloor y\rfloor-h$ pairs of integers at distance $h$). A smaller second factorial moment gives
exactly the observed under-dispersion. The convergence to the Poisson distribution is therefore extremely slow, and the
evidence for Hypothesis~$U(k)$ lies instead in the agreement with the singular series in
\S\ref{subsec:numerics-correlations}--\S\ref{subsec:numerics-conv}.

\FloatBarrier

\bibliographystyle{alphaq}
\bibliography{references}

\end{document}